\documentclass[12pt]{amsart}
\usepackage{amsmath, amssymb, latexsym, amsthm}
\usepackage{mathrsfs}
\usepackage[all]{xy}
\input epsf

\newcommand{\ed}{

\end{document}
}

      \newenvironment{changemargin}[2]{\begin{list}{}{
         \setlength{\topsep}{0pt}\setlength{\leftmargin}{0pt}
         \setlength{\rightmargin}{0pt}
         \setlength{\listparindent}{\parindent}
         \setlength{\itemindent}{\parindent}
         \setlength{\parsep}{0pt plus 1pt}
         \addtolength{\leftmargin}{#1}\addtolength{\rightmargin}{#2}
         }\item }{\end{list}}

\newcommand{\nInc}{\compactN^{\uparrow n}}
\newcommand{\intvl}[2]{{[#1(#2),#1(#2\!+\!1))}}
\newcommand{\Bdd}{\mathbf{B}}

\newcommand{\grbl}{{\mbox{\textit{\tiny gp}}}}

\newcommand{\compactN}{\cl{\mathbb{N}}} 
\newcommand{\blocks}[2]{\op{cl}_{#2}(#1)}
\newcommand{\blocksplus}[2]{\op{cl}^+_{#2}(#1)}

\newcommand{\bq}{\begin{quote}}
\newcommand{\eq}{\end{quote}}
\newcommand{\cl}[1]{\overline{#1}}
\newcommand{\CH}{the Continuum Hypothesis}

\newcommand{\inv}{^{-1}}
\newcommand{\Cantor}{{\{0,1\}^\N}}

\newcommand{\bof}{\op{\fb}}
\newcommand{\bofF}{\bof(\cF)}
\newcommand{\sr}[2]{{#1}} 
\newcommand{\N}{\mathbb{N}}
\newcommand{\NN}{{\N^{\N}}}
\newcommand{\NNup}{{\N^{\uparrow\N}}}
\newcommand{\PN}{{P(\N)}}
\newcommand{\roth}{{[\N]^{\aleph_0}}}
\newcommand{\Fin}{{[\N]^{<\aleph_0}}}
\newcommand{\ici}{{[\N]^{(\aleph_0,\aleph_0)}}}
\newcommand{\Inc}{{\compactN^{\uparrow\N}}}
\newcommand{\powInc}[1]{{\big(\Inc\big)^{#1}}}

\newcommand{\NcompactN}{{\compactN^\N}}
\newcommand{\seq}[1]{\{#1\}_{n\in\N}}
\newcommand{\setseq}[1]{\{#1 : n\in\N\}}
\newcommand{\fbx}[1]{\fbox{$#1$}}

\newcommand{\op}{\operatorname}
\newcommand{\im}{\op{im}}

\newcommand{\maxfin}{\op{maxfin}}

\newcommand{\cI}{\mathcal{I}}
\newcommand{\cJ}{\mathcal{J}}
\newcommand{\scrA}{\mathscr{A}}
\newcommand{\scrB}{\mathscr{B}}
\newcommand{\scrC}{\mathscr{C}}

\newcommand{\B}{\mathcal{B}}
\newcommand{\cB}{\mathcal{B}}

\newcommand{\Tau}{\mathrm{T}}
\newcommand{\cA}{\mathcal{A}}
\newcommand{\cF}{\mathcal{F}}
\newcommand{\cS}{\mathcal{S}}
\newcommand{\cG}{\mathcal{G}}

\newcommand{\cM}{\mathcal{M}}

\newcommand{\cO}{\mathcal{O}}

\newcommand{\Q}{\mathbb{Q}}
\newcommand{\R}{\mathbb{R}}
\newcommand{\cU}{\mathcal{U}}
\newcommand{\Union}{\bigcup}
\newcommand{\cV}{\mathcal{V}}

\newcommand{\Impl}{\Rightarrow}
\long\def\forget#1\forgotten{}

\newcommand{\fb}{\mathfrak{b}}
\newcommand{\ft}{\mathfrak{t}}

\newcommand{\fc}{\mathfrak{c}}
\newcommand{\fd}{\mathfrak{d}}
\newcommand{\fg}{\mathfrak{g}}

\newcommand{\fu}{\mathfrak{u}}

\newcommand{\p}{\mathfrak{p}}

\newcommand{\s}{\mathfrak{s}}

\newcommand{\x}{\times}

\newcommand\comp{^{\text{\tt c}}}
\newcommand{\nin}{\not\in}

\newcommand{\sbst}{\subseteq}
\newcommand{\spst}{\supseteq}
\newcommand{\sm}{\setminus}
\newcommand{\as}{\subseteq^*}
\newcommand{\rest}{\restriction}

\newcommand{\cov}{\op{cov}}
\newcommand{\add}{\op{add}}

\newcommand{\cf}{\op{cf}}
\newcommand{\non}{\op{non}}

\newtheorem{thm}{Theorem}[section]
\newtheorem{prop}[thm]{Proposition}

\newtheorem{lem}[thm]{Lemma}
\newtheorem{cor}[thm]{Corollary}

\theoremstyle{definition}
\newtheorem{defn}[thm]{Definition}
\theoremstyle{remark}
\newtheorem{rem}[thm]{Remark}
\newcommand{\be}{\begin{enumerate}}
\newcommand{\ee}{\end{enumerate}}
\newcommand{\bi}{\begin{itemize}}
\newcommand{\ei}{\end{itemize}}
\newcommand{\bdesc}{\begin{description}}
\newcommand{\edesc}{\end{description}}

\forget
\forgotten

\newcommand{\sone}{\mathsf{S}_1}
\newcommand{\sfin}{\mathsf{S}_{fin}}
\newcommand{\ufin}{\mathsf{U}_{fin}}

\author{Boaz Tsaban}
\thanks{Patially supported by the Koshland Center for Basic Research.}
\address[Boaz Tsaban]{Department of Mathematics, Bar-Ilan University,
Ramat-Gan 52900, Israel;
and
Department of Mathematics,
Weizmann Institute of Science, Rehovot 76100, Israel}
\email{tsaban@math.biu.ac.il}
\urladdr{http://www.cs.biu.ac.il/\~{}tsaban}

\author{Lyubomyr Zdomskyy}
\address[Lyubomyr Zdomskyy]{Department of Mechanics and Mathematics,
Ivan Franko Lviv National University,
Universytetska 1, Lviv 79000, Ukraine; and
Department of Mathematics,
Weizmann Institute of Science,
Rehovot 76100, Israel.}
\curraddr{Kurt G\"odel Research Center for Mathematical Logic, W\"ahringer Str.\ 25, A-1090 Vienna, Austria.}
\email{lzdomsky@gmail.com}

\title[On a problem of Hurewicz]{Scales, fields, and a problem of Hurewicz}

\begin{document}

\begin{abstract}
Menger's basis property is a generalization of
$\sigma$-comp\-actness and admits an elegant combinatorial
interpretation. We introduce a general combinatorial method to
construct non $\sigma$-compact sets of reals with Menger's
property. Special instances of these constructions give known
counterexamples to conjectures of Menger and Hurewicz. We obtain
the first explicit solution to the Hurewicz 1927 problem,
that was previously solved by Chaber and Pol on a dichotomic
basis.

The constructed sets generate nontrivial subfields
of the real line with strong combinatorial properties, and
most of our results can be stated in a Ramsey-theoretic manner.

Since we believe that this paper is of interest to a diverse mathematical
audience, we have made a special effort to make it self-contained and accessible.
\end{abstract}


\maketitle

\emph{Whenever you can settle a question by explicit construction, be not satisfied with purely
existential arguments.}

\medskip

\rightline{Hermann Weyl, Princeton Conference 1946}

\bigskip


\section{Introduction and summary}

Menger's property (1924) is a generalization of
$\sigma$-compactness. Men\-ger conjectured that his property
actually characterizes $\sigma$-compactness. Hurewicz found an
elegant combinatorial interpretation of Menger's property, and
introduced a formally stronger property (1925, 1927). Hurewicz's
property is also implied by $\sigma$-compactness, and Hurewicz
conjectured that his formally stronger property characterizes
$\sigma$-compa\-ctness. He posed the question whether his property
is \emph{strictly} stronger than Menger's. We will call this
question the \emph{Hurewicz Problem}.

In Section \ref{Menger} we define the Menger and Hurewicz properties,
and show that they are extremal
cases of a large family of properties.
We treat this family in a unified manner and obtain,
using a combinatorial approach, many counterexamples to
the above mentioned Conjectures of Menger and Hurewicz.

In Section \ref{HurewiczProblem} we show that a theorem of
Chaber and Pol implies a positive solution to the Hurewicz
problem.
In fact, it establishes the existence of
a set of reals $X$ without the Hurewicz property,
such that all finite powers of $X$ have the Menger property.
However, this solution does not point out a concrete example.
We construct a \emph{concrete} set having Menger's but not Hurewicz's property,
yielding a more elegant and direct solution.

Chaber and Pol's proof is topological.
In Section \ref{FinPows} we show how to obtain Chaber and Pol's result and
extensions of it using the combinatorial approach. In Section \ref{algebra}
we use these results to generate fields (in the algebraic
sense) which are counterexamples to the Hurewicz and Menger
Conjectures and examples for the Hurewicz Problem.
In Section \ref{UM} it is shown that some of our examples are very small, both in the
sense of measure and in the sense of category.

Section \ref{SP} reveals the underlying connections
with the field of selection principles,
where our main results are extended further.
In Section \ref{SFA} we explain how to extend some of the results further,
and in Section \ref{Ramsey} we translate our results into the language of Ramsey theory,
and indicate an application to the undecidable notion of strong measure zero.

\section{The Menger property}\label{Menger}

\subsection{Menger's property and bounded images}\label{menger}
In 1924 Menger introduced the following
basis property for a metric space $X$ \cite{MENGER}:
\begin{quote}
For each basis $\B$ of $X$, there exists a sequence $\seq{B_n}$ in
$\B$ such that $\lim_{n\to\infty}\op{diam}(B_n) = 0$ and $X=\Union_nB_n$.
\end{quote}
Each $\sigma$-compact metric space has this property, and
Menger conjectured that this property characterizes $\sigma$-compactness.
The task of settling this conjecture without special hypotheses was first achieved
in Fremlin and Miller's 1988 paper \cite{FM}, alas in an existential manner.
Concrete counterexamples were given much later \cite{ideals}.
In Section \ref{BFsets}, we describe a general method to produce counterexamples to this conjecture.

In 1927 Hurewicz obtained the following characterization of Menger's property.
Let $\N$ denote the (discrete) space of natural numbers, including $0$,
and endow the \emph{Baire space} $\NN$ with the Tychonoff product topology.
Define a partial order\footnote{
By \emph{partial order} we mean a reflexive and transitive relation.
We do not require its being antisymmetric.}
$\le^*$ on $\NN$ by:
$$f\le^* g\quad \mbox{if}\quad f(n)\le g(n)\mbox{ for all but finitely many }n.$$
A subset $D$ of $\NN$ is \emph{dominating} if
for each $g\in\NN$ there exists $f\in D$ such that $g\le^* f$.
\begin{thm}[Hurewicz \cite{HURE27}]\label{hure}
A set of reals $X$ has Menger's property if, and only if, no
continuous image of $X$ in $\NN$ is dominating.
\end{thm}

Menger's property is a specific instance of a general scheme of properties.

\begin{defn}
For $A,B\sbst\N$, $A\as B$ means that $A\sm B$ is finite.
Let $\roth$ denote the collection of all infinite sets of natural numbers.
A nonempty family $\cF\sbst\roth$ is a \emph{semifilter} if
for each $A\in\cF$ and each $B\sbst\N$ such that $A\as B$, $B\in\cF$ too.
(Note that all elements of $\cF$ are infinite, and $\cF$ is closed under finite modifications
of its elements.)
$\cF$ is a \emph{filter} if it is a semifilter and it is closed under
finite intersections (this is often called a \emph{free} filter).
For $\cF\sbst\roth$ and $f,g\in\NN$, define:
\begin{eqnarray*}
[f\le g] & = & \{n : f(n)\le g(n)\};\\
f\le_\cF g & \mbox{if} & [f\le g]\in\cF.
\end{eqnarray*}
Fix a semifilter $\cF$. A set of reals $X$ satisfies $\Bdd(\cF)$
if each continuous image of $X$ in $\NN$ is bounded with respect to $\le_\cF$,
that is, there is $g\in\NN$ such that for each $f$ in the image of $X$,
$f\le_\cF g$.
\end{defn}

Thus, Menger's property is the same as $\Bdd(\roth)$, and it is the weakest among
the properties $\Bdd(\cF)$ where $\cF$ is a semifilter.

Hurewicz also considered the following property (the \emph{Hurewicz property}) \cite{HURE25}:
Each continuous image of $X$ in $\NN$ is bounded with respect to $\le^*$.
This is also a special case of $\Bdd(\cF)$, obtained when
$\cF$ is the \emph{Fr\'echet filter} consisting of all cofinite sets
of natural numbers. The Hurewicz property is the strongest among
the properties $\Bdd(\cF)$ where $\cF$ is a semifilter, and Hurewicz conjectured
that it characterizes $\sigma$-compactness. This was first disproved by
Just, Miller, Scheepers and Szeptycki in \cite{coc2}, and will also follow from the results below.

The following is easy to verify.

\begin{lem}
For each semifilter $\cF$, $\Bdd(\cF)$ is preserved by continuous images
and is hereditary for closed subsets.\hfill\qed
\end{lem}

This allows us to work in any separable, zero-dimensional metric space
instead of working in $\R$. For brevity, we will refer to any space of
this kind as a \emph{set of reals}.
We consider several canonical spaces which carry a
convenient combinatorial structure.

\subsection{The many faces of the Baire space and the Cantor space}\label{guises}
The Baire space $\NN$ and the \emph{Cantor space} $\Cantor$ are equipped with the product topology.
These spaces will appear under various guises in this paper,
in accordance to the required combinatorial structure.
$\PN$, the collection of all subsets of $\N$, is identified with
$\Cantor$ via characteristic functions, and inherits its topology
(so that by definition $\PN$ and $\Cantor$ are homeomorphic).
$\roth$ is a subspace of $\PN$ and is homeomorphic to $\NN$.
In turn, $\roth$ is homeomorphic to its subspace $\ici$ consisting
of the infinite coinfinite sets of natural numbers.
Similarly,
$\NNup$, the collection of all \emph{increasing} elements of $\NN$,
is homeomorphic to $\NN$.

The following compactification of $\NNup$ appears in
\cite{BaShCon2000}: Let $\compactN=\mathbb{N}\cup\{\infty\}$ be
the one-point compactification of $\N$. Let $\Inc$ be the
collection of all nondecreasing elements $f$ of $\NcompactN$ such
that $f(n)<f(n+1)$ whenever $f(n)<\infty$. For each nondecreasing
finite sequence $s$ of natural numbers, define $q_{s} \in \Inc$ by
$q_s(n) = s(n)$ if $n <|s|$, and $q_s(n) = \infty$ otherwise. Let
$Q$ be the collection of all these elements $q_s$. Then $Q$ is
dense in $\Inc=Q\cup\NNup$. $\Inc$ is another guise of the Cantor
space. Let $\Fin$ denote the finite subsets of $\N$.

\begin{lem}\label{homeo}
Define $\Psi:\Inc\to \PN$ by
$$\Psi(f)=\begin{cases}
\im(f) & f\in\NNup\\
\im(s) & f=q_s\in Q
\end{cases}$$
(in short, $\Psi(f)=\im(f)\sm\{\infty\}$).
Then $\Psi$ is a homeomorphism mapping $Q$ onto $\Fin$ and $\NNup$ onto $\roth$.
\hfill\qed
\end{lem}

Lemma \ref{homeo} says that we can identify sets of natural numbers with their increasing
enumerations, and obtain $\Inc$ (where a finite increasing sequence $s$ is identified with $q_s$).
This identification will be used throughout the paper.
When using it, we will denote elements of $\roth$ by lowercase letters to indicate that
we are also treating them as increasing functions. (Otherwise, we use uppercase letters.)
E.g., for $a,b\in\roth$, $a\le_\cF b$ is an assertion concerning the increasing
enumerations of $a$ and $b$. Similarly for $\le^*$, $\not\le^*$, etc.
Also, $\min\{a,b\}$ denotes the function $f(n)=\min\{a(n),b(n)\}$ for each $n$, and similarly
for $\max\{a,b\}$, etc.

We will need the following lemma from \cite{ideals}.
For the reader's convenience, we reproduce its proof.
\begin{lem}\label{powerlemma}
Assume that $Q^k\sbst X^k\sbst\powInc{k}$, and $\Psi:X^k\to\NN$ is continuous on $Q^k$.
Then there exists $g\in\NN$ such that for
all $x_1,\dots,x_k\in X$,
$$[g<\min\{x_1,\dots,x_k\}]\ \sbst\ [\Psi(x_1,\dots,x_k)\le g].$$
\end{lem}
\begin{proof}
For each $A\sbst\Inc$, let $A\rest n=\{x\rest n : x\in A\}$.
For each $n$, let $\nInc=\Inc\rest n$.
For $\sigma\in\nInc$, write $q_\sigma$ for
$q_{\sigma\rest  m}$ where $m=1+\max\{i<n : \sigma(i)<\infty\}$.

If $\sigma\in\nInc$ and $I$ is a
basic open neighborhood of $q_\sigma$, then there exists a natural number $N$
such that for each $x\in\Inc$ with $x\rest n\in I\rest n$ and
$x(n)>N$, $x\in I$.

Fix $n$.
Use the continuity of $\Psi$ on $Q^k$ to choose,
for each $\vec\sigma = (\sigma_1,\dots,\sigma_k) \in \big(\nInc\big)^k$,
a basic open neighborhood
$$I_{\vec\sigma}= I_{\sigma_1}\x \dots\x I_{\sigma_k}\sbst\powInc{k}$$
of $q_{\vec\sigma}=(q_{\sigma_1},\dots,q_{\sigma_k})$ such that
for all $(x_1,\dots,x_k) \in I_{\vec\sigma}\cap X^k$,
$\Psi(x_1,\dots,\allowbreak x_k)\allowbreak(n)=\Psi(q_{\vec\sigma})(n)$.
For each $i=1,\dots,k$, choose $N_i$ such that for all $x\in\Inc$ with $x\rest n\in
I_{\sigma_i}\rest n$ and $x(n)>N_i$, $x\in I_{\sigma_i}$. Define
$N(\vec\sigma)=\max\{N_1,\dots,N_k\}$.

The set $I_{\vec\sigma}^{(n)}=\{(x_1\rest  n,\dots, x_k\rest  n) :
(x_1,\dots,x_k)\in I_{\vec\sigma}\}$ is open in
$\big(\nInc\big)^k$ and the family $\{I_{\vec\sigma}^{(n)}:
\vec\sigma\in \big(\nInc\big)^k\}$ is a cover of
the compact space $\big(\nInc\big)^k$. Take a finite
subcover $\{I_{\vec\sigma_1}^{(n)},\dots,I_{\vec\sigma_m}^{(n)}\}$
of $\big(\nInc\big)^k$. Let $N=\max\{N(\vec
\sigma_1),\dots,N(\vec \sigma_m)\}$, and define
$$g(n)=\max\{N,\Psi(q_{\vec \sigma_1})(n), \dots, \Psi(q_{\vec \sigma_m})(n)\}.$$
For all $x_1,\dots,x_k\in X$, let $i$ be such that
$(x_1\rest n,\dots,x_k\rest n)\in I_{\vec \sigma_i}^{(n)}$. If
$x_1(n),\dots, x_k(n)>N$, then
$\Psi(x_1,\dots,x_k)(n)=\Psi(q_{\vec\sigma_i})(n)\le g(n)$.
\end{proof}

\subsection{Sets of reals satisfying $\Bdd(\cF)$}\label{BFsets}

\begin{defn}\label{bF}
For a semifilter $\cF$,
let $\bofF$ denote the minimal cardinality of a family $Y\sbst\NN$
which is unbounded with respect to $\le_\cF$.
\end{defn}

The most well known instances of Definition \ref{bF} are
$\fd = \bof(\roth)$ (the minimal cardinality of a dominating family),
and $\fb = \bofF$ where $\cF$ is the Fr\'echet filter
(the minimal cardinality of an unbounded family with respect to $\le^*$).
For a collection (or property) $\cI$ of sets of reals, the
\emph{critical cardinality} of $\cI$ is:
$$\non(\cI)=\min\{|X| : X\sbst\R\mbox{ and }X\nin\cI\}.$$

\begin{lem}
For each semifilter $\cF$, $\non(\Bdd(\cF))=\bofF$.\hfill\qed
\end{lem}

The following notion is our basic building block.

\begin{defn}
$S=\{f_\alpha : \alpha<\bofF\}$ is a \emph{$\bofF$-scale} if $S\sbst\NNup$,
$S$ is unbounded with respect to $\le_\cF$, and for each $\alpha<\beta<\bofF$,
$f_\alpha \le_\cF f_\beta$.
\end{defn}

\begin{lem}\label{EScale}
For each semifilter $\cF$, there exists a $\bofF$-scale.
\end{lem}
\begin{proof}
Let $B=\{b_\alpha : \alpha<\bofF\}\sbst\NN$ be unbounded with respect to $\le_\cF$.
By induction on $\alpha<\bofF$, let $g$
be a witness that $\{f_\beta : \beta<\alpha\}$ is bounded with respect to $\le_\cF$,
and take $f_\alpha=\max\{b_\alpha,g\}$.
Then $S=\{f_\alpha : \alpha<\bofF\}$ is a $\bofF$-scale.
\end{proof}

For a semifilter $\cF$, define
$$\cF^+=\{A\sbst\N : A\comp\nin\cF\}.$$
Note that $\cF^{++}=\cF$.

\begin{rem}
Let $R$ be a binary relation on a set $P$.
A subset $S$ of $P$ is \emph{cofinal} with respect to $R$
if for each $p\in P$ there is $s\in S$ such that $p R s$.
A transfinite sequence $\{p_\alpha : \alpha<\kappa\}$ in $P$ is \emph{nondecreasing}
with respect to $R$ if $p_\alpha R p_\beta $ for all $\alpha \le \beta$.
\end{rem}

Recall that a set of reals $X$ is \emph{meager}
(has Baire \emph{first category}) if it is a countable
union of nowhere dense sets.
Since the autohomeomorphism of $\PN$ defined by
$A\mapsto A\comp$ carries $\cF^+$ to $\cF\comp=P(\N)\sm\cF$, we have
that $\cF$ is meager if, and only if, $\cF^+$ is comeager.

\begin{lem}\label{SFlemma}
Assume that $\cF$ is a semifilter. Then:
$A\in\cF^+$ if, and only if, $A\cap B$ is infinite for each $B\in\cF$.
\end{lem}
\begin{proof}
$(\Leftarrow)$ Assume that $A\cap B$ is infinite for each $B\in\cF$.
Since $A\cap A\comp =\emptyset$, necessarily $A\comp\nin\cF$.

$(\Impl)$ If $B\in\cF$ and $A\cap B$ is finite, then $B\as A\comp$;
thus $A\comp\in\cF$.
\end{proof}

\begin{defn}
For a semifilter $\cF$ and $A\in\cF^+$, define
\begin{eqnarray*}
\cF\rest A & = & \{B\cap A : B\in\cF\};\\
\cF_A & = & \{C\sbst\N : (\exists B\in\cF)\ B\cap A\sbst C\}.
\end{eqnarray*}
\end{defn}

\begin{lem}
For each semifilter $\cF$ and each $A\in\cF^+$,
\be
\item $\cF_A$ is the smallest semifilter containing $\cF\rest A$.
\item $\cF\sbst\cF_A$, and if $\cF$ is a filter, then $\cF_A\sbst\cF^+$.
\hfill\qed
\ee
\end{lem}

\begin{thm}\label{SFscale}
Assume that $\cF$ is a semifilter, and $S=\{f_\alpha : \alpha<\bofF\}$ is a $\bofF$-scale.
Let $X=S\cup Q$. Then: For each continuous $\Psi: X\to\NN$,
there exists $A\in\cF^+$ such that $\Psi[X]$ is bounded with respect to $\le_{\cF_A}$.
\end{thm}
\begin{proof}
Let $g\in\NN$ be as in Lemma \ref{powerlemma}.
Since $S$ is unbounded with respect to $\le_\cF$, there exists $\alpha<\bofF$
such that $f_\alpha\not\le_\cF g$, that is, $A:=[g<f_\alpha]\in\cF^+$.
For each $\beta\ge\alpha$, $[f_\alpha\le f_\beta]\in\cF$.
By Lemma \ref{powerlemma},
$$[\Psi(f_\beta)\le g]\spst [g<f_\beta]\spst A\cap [f_\alpha\le f_\beta]\in\cF\rest A.$$
Let $Y=\Psi[\{f_\beta:\beta<\alpha\}\cup Q]$.
Since $|Y|<\bofF$, $Y$ is $\le_\cF$-bounded by some $h\in\NN$, and we may require
that $[g\le h]=\N$.
Then for each $x\in X$,
$\Psi(x)\le_{\cF_A} h$.
\end{proof}

\begin{cor}\label{F+}
In the notation of Theorem \ref{SFscale}, if
$\cF$ is a filter, then
$X$ satisfies $\Bdd(\cF^+)$.\hfill\qed
\end{cor}

In many cases (including the classical ones),
Theorem \ref{SFscale} implies the stronger assertion
that $X$ satisfies $\Bdd(\cF)$.

\begin{cor}\label{BF}
In the notation of Theorem \ref{SFscale}, assume that
\be
\item $\cF$ is an ultrafilter, or
\item $\cF=\roth$ (Menger property), or
\item $\cF$ is the Fr\'echet filter (Hurewicz property).
\ee
Then $X$ satisfies $\Bdd(\cF)$.
\end{cor}
\begin{proof}
(1) If $\cF$ is an ultrafilter, then $\cF^+=\cF$, and by Corollary \ref{F+}, $X$ satisfies $\Bdd(\cF)$.

(2) If $\cF=\roth$, then for each $A\in\cF^+$, $A$ is cofinite and therefore $\cF_A=\cF$,
so $X$ satisfies $\Bdd(\cF)$.

(3) If $\cF$ is the Fr\'echet filter, then each continuous
image of $X$ in $\NN$ is $\le^*$-bounded when restricted to the infinite set $A$.
To complete the proof, we make the following easy observations.

\begin{lem}\label{wloginc}
The mapping $\Psi:\NN\to\NNup$ defined by
$\Psi(f)(n)=n+f(0)+\dots+f(n)$ is a homeomorphism and preserves $\le^*$-unboundedness.\hfill\qed
\end{lem}

\begin{lem}
If a subset of $\NNup$ is $\le^*$-bounded when restricted to some infinite set $a\sbst\N$,
then it is $\le^*$-bounded.
\end{lem}
\begin{proof}
If $a\as[f\le g]$ for each $f\in Y$ and $g$ is increasing, then
$f\le^* f\circ a \le^* g\circ a$ for each $f\in Y$.
\end{proof}

It follows that each continuous image of $X$ is $\le^*$-bounded,
so $X$ satisfies $\Bdd(\cF)$.
\end{proof}

Items (2) and (3) in Corollary \ref{BF} were first proved in \cite{ideals},
using two specialized proofs.

None of the examples provided by Theorem \ref{SFscale} is trivial:
Each of them is a counterexample to the Menger Conjecture, and some of them
are counterexamples to the Hurewicz Conjecture (see also Section \ref{manycounter}).
Let $\kappa$ be an infinite cardinal.
A set of reals $X$ is \emph{$\kappa$-concentrated} on a set $Q$ if,
for each open set $U$ containing $Q$, $|X\sm U|<\kappa$.
Recall that a set of reals is \emph{perfect} if it is nonempty, closed, and
has no isolated points.

\begin{lem}\label{concen}
Assume that a set of reals $X$ is $\fc$-concentrated on a countable set $Q$.
Then $X$ does not contain a perfect set.
\end{lem}
\begin{proof}
Assume that $X$ contains a perfect set $P$.
Then $P\sm Q$ is Borel and uncountable,
and thus contains a perfect set $C$.
Then $U=\R\sm C$ is open and contains $Q$, and $C=P\sm U\sbst X\sm U$ has cardinality $\fc$.
Thus, $X$ is not $\fc$-concentrated on $Q$.
\end{proof}

\begin{thm}\label{noperfect}
In the notation of Theorem \ref{SFscale}, $X$ does not contain a perfect subset.
In particular, $X$ is not $\sigma$-compact.
\end{thm}
\begin{proof}
If $U$ is an open set containing $Q$, then $K=\big(\Inc\big)\sm U$
is a closed and therefore compact subset of $\Inc$. Thus, $K$ is a
compact subset of $\NNup$, and therefore it is bounded with
respect to $\le^*$. Thus, $|S\cap K|<\bofF$. This shows that $X$
is $\bofF$-concentrated (in particular, $\fc$-concentrated) on $Q$.
Use Lemma \ref{concen}.
\end{proof}

\begin{rem}
If fact, the proof of Theorem \ref{noperfect} gives more:
Since $\bofF\le\fd$, $X$ is $\fd$-concentrated on $Q$, and therefore
\cite{ideals} $X$ has the property $\sone(\Gamma,\cO)$ defined in \cite{coc2}
(see the forthcoming Section \ref{SP}).
By \cite{coc2}, $\sone(\Gamma,\cO)$ is preserved under continuous images and
implies that there are no perfect subsets.
It follows that no continuous image of $X$ contains a perfect subset.
\end{rem}

\subsection{Cofinal scales}
In some situations the following is useful.

\begin{defn}
For a semifilter $\cF$, say that $S=\{f_\alpha : \alpha<\bofF\}$ is a \emph{cofinal $\bofF$-scale}
if:
\be
\item For all $\alpha<\beta<\bofF$, $f_\alpha \le_\cF f_\beta$;
\item For each $g\in\NN$, there is $\alpha<\bofF$ such that for each $\beta\ge \alpha$,
$g\le_\cF f_\beta$.
\ee
\end{defn}

If $\cF\sbst\cF^+$ (in particular, if $\cF$ is a filter),
then every cofinal $\bofF$-scale is a $\bofF$-scale.
If $\cF^+$ is a filter, then every $\bofF$-scale is a cofinal $\bofF$-scale.
Thus, for ultrafilters the notions coincide.

\begin{lem}
Assume that $\cF$ is a semifilter and $\fb(\cF)=\fd$.
Then there exists a cofinal $\bofF$-scale.
\end{lem}
\begin{proof}
Fix a dominating family $\{d_\alpha : \alpha<\fd\}\sbst\NN$.
At step $\alpha<\fd$, choose $f_\alpha\in\NNup$ which is
an upper bound of $\{f_\beta, d_\beta : \beta<\alpha\}$ with respect to
$\le_\cF$ (this is possible because $\fd=\bofF$).
Take $S=\{f_\alpha : \alpha<\fd\}$.

Let $g\in\NN$. Take $\alpha<\fd$ such that $g\le^* d_\alpha$.
For each $\beta\ge\alpha$,  $g\le^* d_\alpha\le_\cF f_\beta$,
and therefore $g\le_\cF f_\beta$.
\end{proof}

\begin{thm}\label{cofscale1}
Assume that $\cF$ is a semifilter.
Then for each cofinal $\bofF$-scale $S=\{f_\alpha : \alpha<\bofF\}$,
$X=S\cup Q$ satisfies $\Bdd(\cF)$.
\end{thm}
\begin{proof}
Assume that $\Psi:X\to\NN$ is continuous.
Let $g\in\NN$ be as in Lemma \ref{powerlemma}.
Take $\alpha<\bofF$ such that for each $\beta\ge\alpha$,
$g\le_\cF f_\beta$. Then for each $\beta\ge\alpha$,
$\Psi(f_\beta)\le_\cF g$.

The cardinality of $\Psi[\{f_\beta : \beta\le \alpha\}\cup Q]$ is smaller
than $\bofF$, and is therefore bounded with respect to $\le_\cF$, either.
It follows that $\Psi[X]$ is bounded with respect to $\le_\cF$.
\end{proof}

\subsection{Many counterexamples to the Hurewicz Conjecture}\label{manycounter}
Recall that Hurewicz conjectured that for sets of reals,
the Hurewicz property is equivalent to $\sigma$-compactness.
In the previous section we gave one type of counterexample,
derived from a $\bofF$-scale where $\cF$ is the Fr\'echet filter.
We extend this construction to a family of semifilters.

\begin{defn}
A family $\cF\sbst\roth$ is \emph{feeble} if there exists
$h\in\NNup$ such that for each $A\in\cF$, $A\cap \intvl{h}{n}\neq\emptyset$
for all but finitely many $n$.
\end{defn}

By a result of Talagrand (see \cite[5.4.1]{CSF}), a semifilter $\cF$ is feeble if,
and only if, it is a meager subset of $\roth$.

\begin{lem}\label{feeblebdd}
Assume that $\cF\sbst\roth$ is feeble.
If $Y\sbst\NNup$ is bounded with respect to $\le_\cF$,
then $Y$ is bounded with respect to $\le^*$.
\end{lem}
\begin{proof}
Take $h\in\NNup$ witnessing that $\cF$ is feeble,
and $g\in\NNup$ witnessing that $Y$ is bounded with respect to $\le_\cF$.
Define $\tilde g\in\NN$ by
$\tilde g(k)=g(h(n+2))$ for each $k\in \intvl{h}{n}$.
It is easy to see that for each $f\in Y$, $f\le^*\tilde g$.
\end{proof}

\begin{cor}\label{bofFeeble}
Assume that $\cF$ is a feeble semifilter. Then $\bofF=\fb$.\hfill\qed
\end{cor}

\begin{thm}\label{feeblescale}
Assume that $\cF$ is a feeble semifilter, and $S=\{f_\alpha : \alpha<\fb\}$ is a $\bofF$-scale.
Then $X=S\cup Q$ has the Hurewicz property.
\end{thm}
\begin{proof}
This is a careful modification of the proof of Theorem \ref{SFscale}.
Let $h\in\NNup$ witness the feebleness of $\cF$.
Assume that $\Psi: X\to\NN$ is continuous.
We may assume that all elements in $\Psi[X]$ are increasing (see Lemma \ref{wloginc}).
Let $g\in\NNup$ be as in Lemma \ref{powerlemma}.
Define $\tilde g\in\NNup$ by
$\tilde g(k)=g(h(n+2))$ for each $k\in \intvl{h}{n}$.

Since $S$ is unbounded with respect to $\le_\cF$, there exists $\alpha<\fb$
such that $A:=[\tilde g<f_\alpha]\in\cF^+$. In particular, $A$ is infinite.
Let $C = \{n : A\cap [h(n-1),h(n))\neq\emptyset\}$.
For each $\beta\ge\alpha$, $[f_\alpha\le f_\beta]\in\cF$.
For all but finitely many $n\in C$, there are
$m\in [f_\alpha\le f_\beta]\cap \intvl{h}{n}$
and $l\in A\cap [h(n-1),h(n))$,
and therefore
$$g(h(n+1))=\tilde g(l)<f_\alpha(l)\le f_\alpha(m)\le f_\beta(m)\le f_\beta(h(n+1)).$$
In particular, $[g<f_\beta]\cap [h(n+1),h(n+2))\neq\emptyset$.
By Lemma \ref{powerlemma},
$$[\Psi(f_\beta)\le g]\spst [g<f_\beta]\spst^* \{h(n+1) : n\in C\}.$$
Thus, $Y=\{\Psi(f_\beta) : \beta\ge\alpha\}$ is $\le^*$-bounded on an infinite set
and therefore $\le^*$-bounded.

Let $Z=\Psi[\{f_\beta:\beta<\alpha\}\cup Q]$.
Since $|Z|<\fb$, $Z$ is $\le^*$-bounded, and therefore
$\Psi[X] = Y\cup Z$ is $\le^*$-bounded.
\end{proof}

By Theorem \ref{noperfect}, each of the sets $X$ of Theorem \ref{feeblescale} is a counterexample
to the Hurewicz Conjecture.

\subsection{Coherence classes}
We make some order in the large family of properties of the form $\Bdd(\cF)$.

\begin{defn}\label{stsubco}
For $h\in\NNup$ and $A\sbst\N$ let
$$\blocks{A}{h} = \Union\{\intvl{h}{n} : A\cap \intvl{h}{n}\neq\emptyset\}.$$
A semifilter $\cS$ is \emph{strictly subcoherent} to a semifilter $\cF$ if
there exists $h\in\NNup$ such that for each $A\in\cS$, $\blocks{A}{h}\in\cF$
(equivalently, there is a monotone surjection $\varphi:\N\to\N$ such that
$\{\varphi[A] : A\in\cS\}\sbst\{\varphi[A] : A\in\cF\}$).
$\cS$ is \emph{strictly coherent} to $\cF$ if each of them is strictly subcoherent to the other.
\end{defn}

The Fr\'echet filter is strictly subcoherent to any semifilter,
so that a semifilter is feeble exactly when it is strictly coherent to the
Fr\'echet filter.

\begin{lem}\label{comeagerseq}
Each comeager semifilter $\cS$ is strictly coherent to $\roth$.
\end{lem}
\begin{proof}
Clearly, any semifilter is strictly subcoherent to $\roth$.
We prove the other direction.
Since $\cS^+$ is homeomorphic to $\cS\comp$, it is meager and thus feeble.
Let $h\in\NNup$ be a witness for that.
Fix $A\in\roth$ and let $B=\blocks{A}{h}$. Then $B\comp\nin\cS^+$,
and therefore $B\in(\cS^+)^+=\cS$.
\end{proof}

\begin{lem}\label{Phi_h}
Let $h\in\NNup$.
Define a mapping $\Phi_h:\NN\to\NN$ by $\Phi_h(f)=\tilde f$,
where for each $n$ and each $k\in\intvl{h}{n}$,
$$\tilde f(k)=\max\{f(i) : i\in\intvl{h}{n}\}.$$
Then:
\be
\item $\Phi_h$ is continuous.
\item For each $f\in\NN$, if $\tilde f=\Phi_h(f)$, then $[f\le \tilde f]=\N$.
\item For each $f,g\in\NN$, if $\tilde f=\Phi_h(f)$, $\tilde g=\Phi_h(g)$,
and $A=[\tilde f\le\tilde g]$, then $A=\blocks{A}{h}$.\hfill\qed
\ee
\end{lem}

\begin{thm}\label{subco}
Assume that $\cS$ and $\cF$ are semifilters such that $\cS$ is strictly subcoherent to
$\cF$.
Then $\Bdd(\cS)$ implies $\Bdd(\cF)$.
In particular, the properties $\Bdd(\cF)$ depend only on the strict-coherence
class of $\cF$.
\end{thm}
\begin{proof}
Assume that $X$ satisfies $\Bdd(\cS)$, and
let $h\in\NNup$ be a witness for $\cS$ being strictly subcoherent to
$\cF$.
Let $Y$ be any continuous image of $X$ in $\NN$.
Then $Y$ satisfies $\Bdd(\cS)$, and therefore so does
$\tilde Y=\Phi_h[Y]$ (where $\Phi_h$ is as in Lemma \ref{Phi_h}), a continuous image of $Y$.
Let $g\in\NN$ be a witness for that, and
take $\tilde g=\Phi_h(g)$.
For each $f\in Y$,
$$[f\le \tilde g]\spst [\tilde f\le \tilde g]\spst [\tilde f\le g]\in\cS,$$
so taking $A=[\tilde f\le \tilde g]$, we have that $[f\le \tilde g]\spst A=\blocks{A}{h}\in\cF$.
\end{proof}

Note that if $\cS$ is a feeble semifilter, then by Theorem
\ref{subco}, $\Bdd(\cS)=\Bdd(\cF)$ where $\cF$ is the Fr\'echet
filter (thus, each set of reals satisfying $\Bdd(\cS)$ is a
counterexample to the Hurewicz Conjecture). In particular,
$\bof(\cS)=\non(\Bdd(\cS))=\non(\Bdd(\cF))=\fb$. Thus, Theorem
\ref{subco} can be viewed as a structural counterpart of Corollary
\ref{bofFeeble}.

\section{A problem of Hurewicz}\label{HurewiczProblem}

\subsection{History}
In his 1927 paper \cite{HURE27}, Hurewicz writes (page 196, footnote 1):
\begin{quote}
Aus der Eigenschaft $E^{**}$ folgt offenbar die Eigenschaft $E^*$.
Die Frage nach der Existenz von Mengen mit der Eigenschaft $E^*$ ohne Eigenschaft $E^{**}$
bleibt hier offen.
\end{quote}
In our language and terminology this reads:
``The Menger property obviously follows from the Hurewicz property.
The question about the existence of sets with the Menger property
and without the Hurewicz property remains open.''

At the correction stage, Hurewicz added there
that Sierpi\'nski proved that the answer is positive if we assume \CH{}.
Thus, the answer is consistently positive. But it remained open whether the answer is
provably positive.

This problem of Hurewicz also appears twice in Lelek's 1969 paper \cite{Lelek69} (pages 210 and 211),
as well as in several recent accounts, for example:
Problem 3 in Just, Miller, Scheepers and Szeptycki's \cite{coc2}.

An existential solution to the Hurewicz Problem was essentially
established by Chaber and Pol at the end of 2002.

\begin{thm}[Chaber-Pol  \cite{ChaPol}]\label{CPThm}
There exists $X\sbst\NN$
such that all finite powers of $X$ have the Menger property, and
$X$ is not contained in any $\sigma$-compact subset of $\NN$.
\end{thm}

In Theorem 5.7 of \cite{coc2} it is proved that a set of reals $X$ has the Hurewicz
property if, and
only if, for each $G_\delta$ set $G$ containing $X$, there is a $\sigma$-compact set $K$
such that $X\sbst K\sbst G$.
Consequently, Chaber and Pol's result implies a positive answer to the Hurewicz Problem,
even when all finite powers of $X$ are required to have the Menger property.

Prior to the present investigation, it was not observed that
the Chaber-Pol Theorem \ref{CPThm} solves the Hurewicz Problem,
and the Hurewi\-cz Problem continued to be raised, e.g.: Problem 1 in
Bukovsk\'y and Hale\v{s}' \cite{BH03}; Problem 2.1 in Bukovsk\'y's
\cite{BukSurv}; Problem 1 in Bukovsk\'y's \cite{Buk03}; Problem
5.1 in the first author's \cite{futurespm}.

Chaber and Pol's solution is existential in the sense that
their proof does not point out a specific example for a set $X$, but instead
gives one example if $\fb=\fd$ (the interesting case), and another if $\fb<\fd$ (the trivial case).
In the current context, this approach was originated in Fremlin and Miller's \cite{FM},
improved in Just, Miller, Scheepers and Szeptycki's  \cite{coc2}
and exploited further in Chaber and Pol's argument.

We will give an explicit solution to the Hurewicz Problem.

\subsection{A solution of the Hurewicz Problem by direct construction}
A continuous metrizable image of the Baire space $\NN$ is called \emph{analytic}.

\begin{lem}\label{analytic}
Assume that $\cA$ is an analytic subset of $\roth$.
Then the smallest semifilter $\cF$ containing $\cA$ is analytic.
\end{lem}
\begin{proof}
For a finite subset $F$ of $\N$, define $\Phi_F:\roth\to \roth$ by
$\Phi_F(A)=A\sm F$ for each $A\in\roth$.
Then $\Phi_F$ is continuous, and therefore $\Phi_F(\cA)$ is analytic.
Let
$$\cB = \Union_{\mbox{\tiny finite }F\sbst\N}\Phi_F(\cA).$$
Then $\cB$ is analytic, and therefore so is
$\cB\x \PN$. Since the mapping $\Phi: \PN\x \PN\to \PN$ defined
by $(A,B)\mapsto A\cup B$ is continuous, we have that $\cF=\Phi[\cB\x \PN]$
is analytic.
\end{proof}

\begin{lem}\label{nonM1}
Assume that $\cF$ is a nonmeager semifilter,
and $Y\sbst\NN$ is analytic.
If $Y$ is bounded with respect to $\le_{\cF^+}$, then
$Y$ is bounded with respect to $\le^*$.
\end{lem}
\begin{proof}
Let $g$ be a $\le_{\cF^+}$-bound of $Y$.
Define $\Phi:Y\to\roth$, by
$$\Phi(f)= [f\le g].$$
Then $\Phi$ is continuous.
Thus, $\Phi[Y]$ is analytic and by Lemma \ref{analytic},
the smallest semifilter $\cS$ containing it is analytic, too.
Since $\cS$ is closed under finite modifications of its elements,
we have by the Topological 0-1 Law \cite[8.47]{Kechris}
that $\cS$ is either meager or comeager.

Note that $\cS\subset\cF^+$.
Since $\cF$ is not meager, $\cF^+$ is not comeager, hence $\cS$ is meager (and therefore feeble).
As $Y$ is bounded with respect to $\le_\cS$ (as witnessed by $g$),
it follows from Lemma \ref{feeblebdd} that $Y$ is bounded with respect to $\le^*$.
\end{proof}

\begin{cor}\label{nonM}
Assume that $\cF$ is a nonmeager semifilter.
Then for each $f\in\NNup$,
the set $\{g\in\NNup : f\leq_\cF g\}$ is nonmeager.
\end{cor}
\begin{proof}
Assume that $\{g\in\NNup : f\leq_\cF g\}$ is meager.
Then there exists a dense $G_\delta$ set $G\sbst\NNup$
such that $g\leq_{\cF^+}f$ for all $g\in G$.
By Lemma \ref{nonM1}, $G$ is bounded with respect to $\le^*$,
and therefore meager; a contradiction.
\end{proof}

\begin{defn}\label{nonM2}
A semifilter $\cS$ is \emph{nonmeager-bounding} if
for each family $Y\sbst\NN$ with $|Y|<\bof(S)$,
the set $\{g\in\NNup : (\forall f\in Y)\ f\leq_\cS g\}$ is nonmeager.
\end{defn}

We will use the following generalization of Definition \ref{stsubco}.

\begin{defn}
For $h\in\NNup$ and $A\sbst\N$ let
$$\blocksplus{A}{h} = \Union\{[h(n),h(n+3)) : A\cap [h(n+1),h(n+2))\neq\emptyset\}.$$
A semifilter $\cS$ is \emph{subcoherent} to a semifilter $\cF$ if
there exists $h\in\NNup$ such that for each $A\in\cS$, $\blocksplus{A}{h}\in\cF$.
$\cS$ is \emph{coherent} to $\cF$ if each of them is subcoherent to the other.
\end{defn}

It is often, but not always, the case that subcoherence coincides with strict subcoherence---see Chapter 5 of \cite{CSF}.

\begin{prop}\label{nonMbounding}
Assume that $\cS$ is a semifilter. If any of the following holds, then
$\cS$ is nonmeager-bounding:
\be
\item $\cS$ is a nonmeager filter, or
\item $\cS=\roth$, or
\item $\cS$ is coherent to a nonmeager filter, or
\item $\cS$ is comeager.
\ee
\end{prop}
\begin{proof}
(1) Assume that $Y\sbst\NNup$ and $|Y|<\bof(\cS)$. Let $f\in\NNup$ be a $\le_\cS$-bound of $Y$.
By Corollary \ref{nonM}, $\{g\in\NNup : f\leq_\cS g\}$ is nonmeager.
Since $\cS$ is a filter, $\le_\cS$ is transitive, and therefore
each member in this nonmeager set is a $\le_\cS$-bound of $Y$.

(2) Assume that $Y\sbst\NNup$ and $|Y|<\fd$.
We may assume that $Y$ is closed under pointwise maxima.
Let $g\in\NNup$ be a witness for the fact that $Y$ is not dominating.
Then $\{[f\le g] : f\in Y\}$ is closed under taking finite intersections.
Let $\cU$ be an ultrafilter extending it.
By (1), $Z=\{h\in\NNup : g\le_\cU h\}$ is nonmeager.
For each $h\in Z$, $f\le_\cU g\le_\cU h$ (in particular, $f\le_\roth g$)
for each $f\in Y$.

(3) Any semifilter coherent to a \emph{filter} is actually strictly coherent to it \cite[5.5.3]{CSF}.
Thus, assume that $\cS$ is \emph{strictly} coherent to a nonmeager filter $\cF$.
Then there is a monotone surjection $\varphi:\N\to\N$ such that
$\{\varphi[A] : A\in\cS\}=\{\varphi[A] : A\in\cF\}$ \cite[5.5.2]{CSF}.
The filter $\cG$ generated by $\{\varphi\inv[\varphi[A]] : A\in\cS\}$
is contained in $\cS$. Since $\cG$ is coherent to $\cS$, it is nonmeager
and $\bof(\cG)=\bof(\cS)$ \cite[5.3.1 and 10.1.13]{CSF}.
Thus, $\cG$ is nonmeager-bounding and since $\bof(\cG)=\bof(\cS)$ and $\le_\cS$ extends $\le_\cG$,
$\cS$ is nonmeager-bounding.

(4) Using Lemma \ref{comeagerseq},
let $h\in\NNup$ be a witness for $\roth$ being strictly subcoherent to $\cS$.
Assume that $Y\sbst\NN$ and $|Y|<\bof(\cS)$.
For each $f\in Y$ define $\tilde f\in\NN$ by
$\tilde f(n)=\max\{f(k) : k\in\intvl{h}{n}\}$.
By (2),
$Z=\{g\in\NNup : (\forall f\in Y)\ \tilde f\le_\roth g\}$ is nonmeager.
Fix any $g$ in this nonmeager set.
Let $f\in Y$, and
$A=[\tilde f\le g]$. $A\in\roth$, and for each
$n\in A$ and each $k\in\intvl{h}{n}$,
$$f(k)\le\tilde f(n)\le g(n)\le g(h(n))\le g(k),$$
that is, $[f\le g]\spst\Union_{n\in A}\intvl{h}{n}$.
By Lemma \ref{comeagerseq}, the last set is a member of $\cS$.
\end{proof}

\begin{rem}
Under some set theoretic hypotheses, e.g., $\fb=\fd$ or $\fu<\fg$, all nonmeager
semifilters are nonmeager-bounding.%
\end{rem}

The assumptions on $\cF$ in the following theorem hold for $\cF=\roth$.
Thus, this theorem implies the promised solution to the Hurewicz Problem.

\begin{thm}\label{answerHure}
Assume that $\cF$ is a nonmeager-bounding semifilter with $\bofF=\fd$.
Then there is a cofinal $\bofF$-scale $S=\{f_\alpha : \alpha<\fd\}$
such that the set $X=S\cup Q$ satisfies $\Bdd(\cF)$ but does not have the Hurewicz property.
\end{thm}
\begin{proof}
We will identify $\Inc$ with $\PN$, identifying $Q$ with $\Fin$ and $\NNup$ with $\roth$
(see Lemma \ref{homeo} and the discussion following it).
Recall that $\ici$ is the collection of infinite coinfinite subsets of $\N$.
For each $g\in\NN$, $\{a\in\ici : a\le^* g\}$ is meager,
and therefore so is $M_g:=\{a\in\ici : a\comp\le^* g\}$ (since $A\mapsto A\comp$ is
an autohomeomorphism of $\ici$).

Fix a dominating family $\{d_\alpha : \alpha<\fd\}\sbst\NN$.
Define $a_\alpha\in\ici$ by induction on $\alpha<\fd$, as follows:
At step $\alpha$ use the fact that $\cF$ is nonmeager-bounding
to find $a_\alpha\in\ici\sm M_{d_\alpha}$ which is
a bound for $\{d_\beta,a_\beta : \beta<\alpha\}$ with respect to $\le_\cF$.
Take $S=\{a_\alpha : \alpha<\fd\}$.

By Theorem \ref{cofscale1}, $X=S\cup Q$ satisfies $\Bdd(\cF)$.
But $\{x\comp : x\in X\}$ is a homeomorphic image of $X$ in $\NN$,
and is unbounded (with respect to $\le^*$),
since for each $\alpha<\fd$, $a_\alpha\comp\not\le^* d_\alpha$.
Thus, $X$ does not have the Hurewicz property.
\end{proof}

The methods that Chaber and Pol used to prove their Theorem \ref{CPThm} are topological.
We proceed to show that Chaber and Pol's Theorem can also be obtained using the combinatorial approach.

\section{Finite powers and the Chaber-Pol Theorem}\label{FinPows}

Having the property $\Bdd(\cF)$ in all finite powers is useful for the generation of
(nontrivial) groups and other algebraic objects satisfying $\Bdd(\cF)$.
In this section we restrict attention to filters.

\begin{thm}\label{Fscale}
Assume that $\cF$ is a filter, and $S=\{f_\alpha : \alpha<\bofF\}$ is a $\bofF$-scale.
Let $X=S\cup Q$. Then: For each $k$ and each continuous $\Psi: X^k\to\NN$,
there exist elements $A_1,\dots,A_k\in\cF^+$ such that $\Psi[X^k]$ is bounded with respect to
$\le_{\cF_{A_1}\cup\dots\cup\cF_{A_k}}$.
\end{thm}
The proof of Theorem \ref{Fscale} is by induction on $k$.
To make the induction step possible, we strengthen its assertion.
\begin{prop}\label{FscaleProp}
For each $k$ and each family $\scrC$ of less than $\bofF$ continuous
functions from $X^k$ to $\NN$,
there exist elements $A_1,\dots,\allowbreak A_k\in\cF^+$ such that $\Union\{\Psi[X^k] : \Psi\in\scrC\}$
is bounded with respect to $\le_{\cF_{A_1}\cup\dots\cup\cF_{A_k}}$.
\end{prop}
\begin{proof}
For each $\Psi\in\scrC$, let $g_\Psi\in\NN$ be as in Lemma \ref{powerlemma}.
Since $|\scrC|<\bofF$, there is $g_0\in\NN$ such that $g_\Psi\le_\cF g_0$ for each $\Psi\in\scrC$.
Choose $\alpha<\bofF$ such that $[g_0<f_\alpha]\in\cF^+$.
Then $A_k:=[g_0<f_\alpha]\in\cF^+$.
We continue by induction on $k$.

$k=1$: By Lemma \ref{powerlemma}, for each $\beta\ge\alpha$ and each $\Psi\in\scrC$,
$[\Psi(f_\beta)\le g_0]\in\cF_{A_1}$.
Since the cardinality of the set
$$\{\Psi(f) : \Psi\in\scrC, f\in \{f_\beta : \beta<\alpha\}\cup Q\}$$
is smaller than $\bofF$, this set is bounded with respect to $\le_\cF$,
by some function $h\in\NN$. Take $g=\max\{g_0,h\}$.

$k=m+1$:
For all $\alpha_1,\dots,\alpha_k\ge\alpha$,
we have by Lemma \ref{powerlemma} that
\begin{eqnarray*}
\lefteqn{[\Psi(f_{\alpha_1},\dots,f_{\alpha_k})\le g_0]\spst}\\
& \spst & [g_0<\min\{f_{\alpha_1},\dots,f_{\alpha_k}\}]\spst
A_k\cap\bigcap_{i=1}^{k}[f_{\alpha}\le f_{\alpha_i}]\in\cF\rest A_k.
\end{eqnarray*}

For each $f\in\{f_\beta : \beta<\alpha\}\cup Q$
and each $i=1,\dots,k$ define
$\Psi_{i,f}:X^m\to\NN$ by
$$\Psi_{i,f}(x_1,\dots,x_m) = \Psi(x_1,\dots,x_{i-1},f,x_i,\dots,x_m).$$
Since there are less than $\bofF$ such functions,
we have by the induction hypothesis $A_1,\dots,A_m\in\cF^+$ such that
$$\Union\{\Psi_{i,f}[X^m] : i=1,\dots,k,\ f\in\{f_\beta : \beta<\alpha\}\cup Q,\ \Psi\in\scrC\}$$
is bounded with respect to $\le_{\cF_{A_1}\cup\dots\cup\cF_{A_m}}$.
Let $h\in\NN$ be such a bound, and take $g=\max\{g_0,h\}$.
Then $\Union\{\Psi[X^k] : \Psi\in\scrC\}$
is bounded with respect to $\le_{\cF_{A_1}\cup\dots\cup\cF_{A_k}}$.
\end{proof}

This completes the proof of Theorem \ref{Fscale}.\hfill\qed

\begin{cor}\label{F+powers}
In the notation of Theorem \ref{Fscale},
all finite powers of $X$ satisfy $\Bdd(\cF^+)$.\hfill\qed
\end{cor}

Item (2) in Corollary \ref{BFpowers} was first proved in \cite{ideals},
using a specialized proof.

\begin{cor}\label{BFpowers}
In the notation of Theorem \ref{Fscale}, assume that
\be
\item $\cF$ is an ultrafilter, or
\item $\cF$ is the Fr\'echet filter (Hurewicz property).
\ee
Then all finite powers of $X$ satisfy $\Bdd(\cF)$.
\end{cor}
\begin{proof}
(1) If $\cF$ is an ultrafilter, then $\cF^+=\cF$.

(2) Fix $k$ and a continuous $\Psi: X^k\to\NN$.
We may assume that $\Psi[X^k]\sbst\NNup$.
Apply Theorem \ref{Fscale}, and let $g\in\NN$ be a witness for
$\Psi[X^k]$ being bounded with respect to $\le_{\cF_{A_1}\cup\dots\cup\cF_{A_k}}$.
For each $i=1,\dots,k$ let $Y_i = \{f\in\Psi[X^k] : f\le_{\cF_{A_i}} g\}$.
Then each $Y_i$ is bounded, and therefore so is
$$\Union_{i=1}^k Y_i = \Psi[X^k].\qedhere$$
\end{proof}

For later use, we point out the following.

\begin{thm}\label{cofscale2}
Assume that $\cF$ is a filter and $S=\{f_\alpha : \alpha<\bofF\}$ is a cofinal $\bofF$-scale.
Then all finite powers of the set $X=S\cup Q$ satisfy $\Bdd(\cF)$.
\end{thm}
\begin{proof}
This is a part of the proof of Theorem \ref{Fscale},
replacing each $\cF^+$ with $\cF$ (in this case the proof can be simplified).
\end{proof}

The cardinal $\fd$ is not provably regular. However, in most of the known models
of set theory it is regular.
In Theorem 16 of \cite{ideals}, a weaker version of Theorem \ref{UFscale} is established
using various hypotheses, all of which imply that $\fd$ is regular.

\begin{thm}\label{UFscale}
Assume that $\fd$ is regular. Then there is an ultrafilter $\cU$
with $\bof(\cU)=\fd$, and a $\bof(\cU)$-scale $S=\{f_\alpha : \alpha<\fd\}$
such that all finite powers of the set $X=S\cup Q$ satisfy $\Bdd(\cU)$,
but $X$ does not have the Hurewicz property.
\end{thm}
\begin{proof}
There always exists an ultrafilter $\cU$ with $\bof(\cU)=\cf(\fd)$ \cite{Canj}.
Since ultrafilters are not meager, we have by
Proposition \ref{nonMbounding}(1) that $\cU$ is nonmeager-bounding.
Take a $\bof(\cU)$-scale $S=\{f_\alpha : \alpha<\fd\}$ as in Theorem
\ref{answerHure}, so that the set $X=S\cup Q$ does not have the
Hurewicz property.
By Corollary \ref{BFpowers}(1), all finite powers of $X$ satisfy $\Bdd(\cU)$.
\end{proof}

\begin{rem}\label{dicho}
In particular, we obtain Chaber and Pol's Theorem \ref{CPThm}:
\be
\item If $\fd$ is regular, use Theorem \ref{UFscale}.
Otherwise, let $X$ be any unbounded subset of $\NN$ of cardinality $\cf(\fd)$.
This proof is still on a dichotomic basis, but the dichotomy here puts more
weight on the interesting case (since $\fb<\cf(\fd)=\fd$ is consistent).
\item The sets in this argument are of cardinality $\cf(\fd)$, while Chaber and Pol's sets are of cardinality $\fb$.
To get sets of cardinality $\fb$, use the dichotomy ``$\fb=\fd$ (which implies
that $\fd$ is regular) or $\fb<\fd$'' instead.
\ee
\end{rem}

\begin{rem}\label{Polish}
Like in Chaber and Pol's \cite{ChaPol},
our constructions can be carried out in any nowhere locally compact Polish space $P$:
Fix a countable dense subset $E$ of $P$. Since $E$ and our $Q$ are both countable metrizable
with no isolated points, they are both homeomorphic to the space $\Q$ of rational numbers,
and hence are homeomorphic via some map $\varphi:Q\to E$. According to
Lavrentiev's Theorem \cite[3.9]{Kechris}, $\varphi$ can be extended to a homeomorphism
between two (dense) $G_\delta$-sets containing $Q$ and $E$, respectively.
Now, every $G_\delta$ set $G$ in $\Inc$ containing $Q$ contains the set $\{f\in\NNup:f\not\le^* g\}$
for some fixed $g\in\NNup$, in which our constructions can be carried out.
\end{rem}

\section{Finite powers for arbitrary feeble semifilters}

We extend Theorem \ref{feeblescale} and Corollary \ref{BFpowers}(2).

\begin{thm}\label{feebleprods}
Assume that $\cF_1,\dots,\cF_k$ are feeble semifilters, and for each $i=1,\dots,k$,
$S_i=\{f^i_\alpha : \alpha<\fb\}$ is a $\bof(\cF_i)$-scale and $X_i=S_i\cup Q$.
Then $\prod_{i=1}^k X_i$ has the Hurewicz property.
\end{thm}
\begin{proof}
The proof is by induction on $k$.
The case $k=1$ is Theorem \ref{feeblescale}, so
assume that the assertion holds for $k-1$ and let us prove it for $k$.

Let $h_1,\dots,h_k\in\NNup$ witness the feebleness of $\cF_1,\dots,\cF_k$.
Take $h\in\NNup$ such that for each $n$ and each $i=1,\dots,k$,
$\intvl{h}{n}$ contains some interval $\intvl{h_i}{j}$.
Clearly, $h$ witnesses the feebleness of all semifilters $\cF_1,\dots,\cF_k$.

Assume that $\Psi: \prod_{i=1}^k X_i\to\NN$ is continuous.
We may assume that all elements in $\Psi[X]$ are increasing.
Let $g\in\NNup$ be as in Lemma \ref{powerlemma}, and
define $\tilde g\in\NNup$ by
$\tilde g(m)=g(h(n+2))$ for each $m\in \intvl{h}{n}$.

\begin{lem}\label{simul}
Assume that $Y_1,\dots,Y_k\sbst\NNup$ are unbounded (with respect to $\le^*$) and $g\in\NN$.
Then there exist $f_i\in Y_i$, $i=1,\dots,k$, such that
$[g<\min\{f_1,\dots,f_k\}]$ is infinite.
\end{lem}
\begin{proof}
Take $f_1\in Y_1$ such that $A_1=[g<f_1]$ is infinite.
As all members of $Y_2$ are increasing and $Y_2$ is unbounded,
$Y_2$ is not bounded on $A_1$, thus there is $f_2\in Y_2$
such that $A_2=[g<\min\{f_1,f_2\}]=A_1\cap[g< f_2]$ is infinite.
Continue inductively.
\end{proof}

Use Lemma \ref{simul} to choose $\alpha_1,\dots,\alpha_k<\fb$ such that
$A=[\tilde g < \min\{f^1_{\alpha_1},\dots,f^k_{\alpha_k}\}]$ is infinite.
Let $C = \{n : A\cap [h(n-1),h(n))\neq\emptyset\}$.
Take $\alpha=\max\{\alpha_1,\dots,\alpha_k\}$.

As in the proof of Theorem \ref{feeblescale}, we have that
for each $\beta\ge\alpha$ and each $i=1,\dots,k$,
$g(h(n+1))<f^i_\beta(h(n+1))$
for all but finitely many $n\in C$.
Thus, for all $\beta_1,\dots,\beta_k\ge\alpha$,
$$g(h(n+1))<\min\{f^1_{\beta_1}(h(n+1)),\dots,f^k_{\beta_k}(h(n+1))\}$$
for all but finitely many $n\in C$.
By Lemma \ref{powerlemma},
$$[\Psi(f^1_{\beta_1},\dots,f^k_{\beta_k})\le g]\spst [g<\min\{f^1_{\beta_1},\dots,f^k_{\beta_k}\}]
\spst^* \{h(n+1) : n\in C\},$$
that is, $\{\Psi(f^1_{\beta_1},\dots,f^k_{\beta_k}) : \beta_1,\dots,\beta_k\ge\alpha\}$ is
$\le^*$-bounded on an infinite set and therefore $\le^*$-bounded.

It follows, as at the end of the proof of Proposition \ref{FscaleProp},
that the image of $\Psi$ is a union of less than $\fb$ many $\le^*$-bounded sets, and
is therefore $\le^*$-bounded.
\end{proof}

\section{Adding an algebraic structure}\label{algebra}

In this section we show that most of our examples can be chosen to have an algebraic
structure.

A classical result of von Neumann \cite{vN} asserts that there exists a subset $C$ of $\R$
which is homeomorphic to the Cantor space and is algebraically independent over $\Q$.
Since the properties $\Bdd(\cF)$ are preserved under continuous images, we may
identify $\Inc$ with such a set $C\sbst\R$, and for $X\sbst\Inc$
consider the subfield $\Q(X)$ of $\R$ generated by $\Q\cup X$.
The following theorem extends Theorem 1 of \cite{o-bdd} significantly.

\begin{thm}\label{basicfields}
Assume that $\cF$ is a filter, $S=\{f_\alpha : \alpha<\bofF\}$ is a $\bofF$-scale, and
$X=S\cup Q$. Then:
\be
\item All finite powers of $\Q(X)$ satisfy $\Bdd(\cF^+)$.
\item If $\cF$ is an ultrafilter, then
all finite powers of $\Q(X)$ satisfy $\Bdd(\cF)$.
\item If $\cF$ is a feeble filter, then
all finite powers of $\Q(X)$ have the Hurewicz property.
\ee
On the other hand,
\be
\item[(4)] For each property $P$ of sets of reals which is hereditary for closed subsets,
if $X$ does not have the property $P$, then $\Q(X)$ does not have the property $P$, either.
\item[(5)] $\Q(X)$ is not $\sigma$-compact.
\ee
\end{thm}
\begin{proof}
(1) Denote by $\Q(t_1,\dots,t_n)$ the field of all rational functions
in the indeterminates $t_1,\dots,t_n$ with coefficients in $\Q$.
For each $n$,
$$\Q_n(X)=\left\{r(x_1,\dots,x_n) : r\in\Q(t_1,\dots,t_n),\ x_1,\dots,x_n\in X\right\}$$
is a union of countably many continuous images of $X^n$, thus for each $k$, $(\Q_n(X))^k$
is a union of countably many continuous images of $X^{nk}$, which by Corollary \ref{F+powers}
satisfy $\Bdd(\cF^+)$.

For a family $\cI$ of sets of reals with $\Union\cI\nin\cI$, let
$$\add(\cI)=\min\{|\cJ| : \cJ\sbst\cI\mbox{ and }\Union\cJ\nin\cI\}.$$

\begin{lem}
For each semifilter $\cF$, $\add(\Bdd(\cF))\ge\fb$.
If $\cF$ is a filter, then $\add(\Bdd(\cF))=\bofF$.\hfill\qed
\end{lem}

Since $\Bdd(\cF^+)$ is preserved under taking continuous images
and countable unions, we have that each set $(\Q_n(X))^k$ satisfies
$\Bdd(\cF^+)$, and therefore so does $(\Q(X))^k=\Union_n(\Q_n(X))^k$.

(2) and (3) are obtained similarly, as consequences of Corollary \ref{BFpowers}
and Theorem \ref{feebleprods}, respectively.

(4) Since $X\sbst C$ and $C$ is algebraically independent, we have that
$\Q(X)\cap C=X$ and therefore $X$ is a closed subset of $\Q(X)$.

(5) Use (4) and apply Theorem \ref{noperfect}.
\end{proof}

The following theorem extends Theorem 5 of \cite{o-bdd}, and shows
that even fields can witness that the Hurewicz property is stronger
than Menger's.

\begin{thm}\label{regdfield}
If $\fd$ is regular, then for the set $X$ of Theorem \ref{UFscale},
all finite powers of $\Q(X)$ have the Menger property, but
$\Q(X)$ does not have the Hurewicz property.
\end{thm}
\begin{proof}
This follows from Theorems \ref{UFscale} and \ref{basicfields}(4).
\end{proof}

The following solves Hurewicz's Problem for subfields of $\R$.
\begin{cor}\label{cfdfield}
There exists $X\sbst\R$ of cardinality $\cf(\fd)$ such that all finite powers of $\Q(X)$
have Menger's property, but $\Q(X)$ does not have the Hurewicz property.
\end{cor}
\begin{proof}
Take the dichotomic examples of Remark \ref{dicho}.
\end{proof}

\begin{rem}
Problem 6 of \cite{o-bdd} and Problem 1.3 of \cite{ict} ask (according to the forthcoming Section \ref{SP})
whether there exists a subgroup $G$
of $\R$ such that $|G|=\fd$ and $G$ has Menger's property.
Theorem \ref{regdfield} answers the question in the affirmative under the additional weak assumption
that $\fd$ is regular. Corollary \ref{cfdfield} answers affirmatively the analogous question where
$\fd$ is replaced by $\cf(\fd)$.
\end{rem}

\begin{rem}
We can make all of our examples subfields of any nondiscrete, separable,
completely metrizable field $\mathbb{F}$.
Examples for such fields are, in addition to $\R$, the complex numbers $\mathbb{C}$, and the $p$-adic numbers $\Q_p$.
More examples involving meromorphic functions or formal Laurent series are available in \cite{Pfef00}.
To this end, we use Mycielski's extension of von Neumann's Theorem, asserting that
for each countable dense subfield $\Q$ of $\mathbb{F}$,
$\mathbb{F}$ contains an algebraically independent (over $\Q$)
homeomorphic copy of the Cantor space (see \cite{Pfef00} for a proof).
\end{rem}

\section{Smallness in the sense of measure and category}\label{UM}

A set of reals $X$ is \emph{null} if it has Lebesgue measure zero.
$X$ is \emph{universally null} if every Borel isomorphic image of $X$ in $\R$
is null. Equivalently, for each finite $\sigma$-additive measure $\mu$
on the Borel subsets of $X$ such that $\mu\{x\}=0$ for each $x\in X$,
$\mu(X)=0$.
A classical result of Marczewski asserts that each product of two universally
null sets of reals is universally null.

A set of reals $X$ is \emph{perfectly meager} if for each perfect set $P$, $X\cap P$ is
meager in the relative topology of $P$.
It is \emph{universally meager} if
each Borel isomorphic image of $X$ in $\R$ is meager.
Zakrzewski \cite{zakrUFC} proved that each product of two universally meager
sets is universally meager.

As in Section \ref{algebra},
we identify the Cantor space with a subset of $\R$ which
is algebraically independent over $\Q$.

\begin{thm}\label{PlewPol}
Let $\cF$ be the Fr\'echet filter,
$S=\{f_\alpha : \alpha<\fb\}$ be any $\bofF$-scale,
and $X=S\cup Q$.
Then:
All finite powers of $\Q(X)$ have the Hurewicz property and are universally null and universally meager.
\end{thm}
\begin{proof}
Theorem \ref{basicfields} deals with the first assertion.

Plewik \cite{Plew93} proved that every set $S$ as above is
both universally null and universally meager.\footnote{The
latter assertion also follows from Corollary \ref{BF} and Theorem \ref{noperfect},
by a result of Zakrzewski \cite{zakrUFC} which asserts that every set of reals having the Hurewicz
property and not containing perfect sets is universally meager.}
Since both properties are preserved under taking countable unions and are
satisfied by singletons, we have that $X$ is universally null and universally meager.
Consequently, all finite powers of $X$ are universally null and universally meager.

We should now understand why these properties would also hold for $\Q(X)$
and its finite powers.
To this end, we use some results of Pfeffer and Prikry.
The presentation is mutatis mutandis the one from Pfeffer's \cite{Pfef00}, in
which full proofs are supplied.

Let
$$\Q'(t_1,\dots,t_n)=\Q(t_1,\dots,t_n)\sm\Union^n_{i=1}\Q(t_1,\dots,t_{i-1},t_{i+1},\dots,t_n).$$
The usual order in the Cantor set induces an order $\preceq$ in $C$,
which is closed in $C^2$. For $X\sbst C$ and
each $m$ and $k$, define
$$X_{m,k}=\{(x_1,\dots, x_m)\in X^m : x_1 \preceq \dots \preceq x_m,\ (\forall i\neq j)\ |x_i-x_j|>1/k\}.$$
Each $X_{m,k}$ is closed in $X^m$,
in particular, each $C_{m,k}$ is compact.
Since $C$ is algebraically independent, each $r\in \Q'(t_1,\dots,t_m)$ defines
a continuous map $(a_1,\dots,a_m)\mapsto r(a_1,\dots,a_m)$ from
$\Union_k C_{m,k}$ to $\R$.
It follows that $r$ is a homeomorphism into $\Q(X)$, and that
$$\Q(X) = \Q\cup \Union_{m,k\in\N}\Union_{r\in \Q'(t_1,\dots,t_m)} r[X_{m,k}].$$

For each $m_1,m_2,k_1,k_2$, $X_{m_1,k_1}\x X_{m_2,k_2}\sbst X^{m_1+m_2}$ and is therefore
universally null and universally meager. Thus, so is each homeomorphic copy
$r_1[X_{m_1,k_1}]\x r_2[X_{m_2,k_2}]$ of $X_{m_1,k_1}\x X_{m_2,k_2}$, where
$r_1\in \Q'(t_1,\dots,t_{m_1}), r_2\in \Q'(t_1,\dots,t_{m_2})$.
A similar assertion holds for products of any finite length.
Consequently, each finite power of $\Q(X)$ is a countable union of
sets which are universally null and universally meager, and is therefore
universally null and universally meager.
\end{proof}

\section{Connections with selection principles}\label{SP}

\subsection{Selection principles}
In his 1925 paper \cite{HURE25}, Hurewicz introduced two
properties of the following type.
For collections $\scrA,\scrB$ of covers of a space $X$, define
\bdesc
\item[$\ufin(\scrA,\scrB)$]
For each sequence $\seq{\cU_n}$ of members of $\scrA$
which do not contain a finite subcover,
there exist finite subsets $\cF_n\sbst\cU_n$, $n\in\N$,
such that $\setseq{\bigcup\cF_n}\in\scrB$.
\edesc

\begin{figure}[!htp]
\begin{center}
\epsfysize=6 truecm {\epsfbox{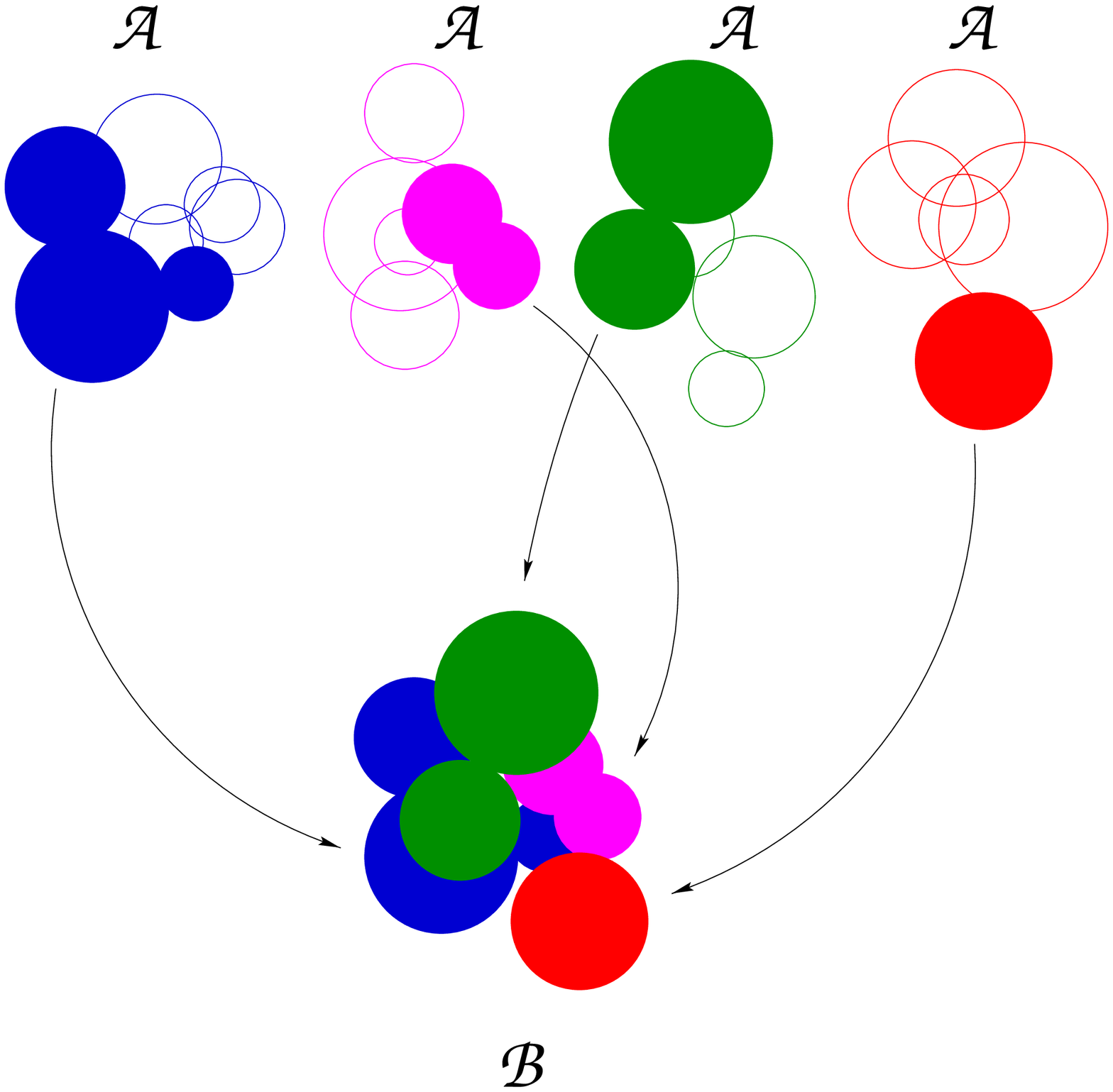}}
\caption{$\ufin(\scrA,\scrB)$}
\end{center}
\end{figure}
Hurewicz (essentially) proved that if $X$ is a set of reals and
$\cO$ is the collection of all open covers
of $X$, then $\ufin(\cO,\cO)$ is equivalent to $\Bdd(\roth)$ (the Menger property).
He also introduced the following property:
Call an open cover $\cU$ of $X$ a \emph{$\gamma$-cover}
if $\cU$ is infinite, and each $x\in X$ belongs to all but finitely many
members of $\cU$.
Let $\Gamma$ denote the collection of all open $\gamma$-covers of $X$.
Hurewicz proved that for sets of reals,
$\ufin(\cO,\Gamma)$ is the same as $\Bdd(\cF)$ where $\cF$ is the
Fr\'echet filter (the Hurewicz property).

Here too, the properties $\ufin(\cO,\cO)$ and $\ufin(\cO,\Gamma)$
are specific instances of a general scheme of properties.
\begin{defn}\mbox{}
\be
\item Let $\cU$ be a cover of $X$ enumerated bijectively
as $\vec\cU=\{U_n : n\in\N\}$.
The \emph{Marczewski characteristic function} of $\vec\cU$,
$h_{\vec\cU}:X\to \PN$, is defined by
$$h_{\vec\cU}(x) = \{ n : x\in U_n\}$$
for each $x\in X$.
\item Let $\cF$ be a semifilter.
\be
\item $\cU$ is an \emph{$\cF$-cover} of $X$ if
there is a bijective enumeration $\vec\cU=\{U_n : n\in\N\}$ such that
$h_{\vec\cU}[X]\sbst\cF$.
\item $\cO_\cF$ is the collection of all open $\cF$-covers of $X$.
\ee
\ee
\end{defn}
If $\cF=\roth$, then it is easy to see that $\ufin(\cO,\cO_\cF)=\ufin(\cO,\cO)$ \cite{coc2}.
If $\cF$ is the Fr\'echet filter, then $\cO_\cF=\Gamma$ and therefore $\ufin(\cO,\cO_\cF)=\ufin(\cO,\Gamma)$.
The families of covers $\cO_\cF$ were first introduced and studied in a similar
context by Garc\'{i}a{}-Ferreira and Tamariz-Mascar\'ua \cite{GT1, GTpseq}.

\begin{defn}
$S_\N$ is the collection of all permutations $\sigma$ of $\N$.
For a semifilter $\cF$ and $\sigma\in S_\N$, write $\sigma\cF=\{\sigma[A] : A\in\cF\}$.
\end{defn}

Observe that if $\cF$ is $\roth$ or the Fr\'echet filter, then
$\sigma\cF=\cF$ for all $\sigma$. Consequently,
the following theorem generalizes Hurewicz's Theorem.

\begin{thm}\label{bddimg}
Assume that $\cF$ is a semifilter. For a set of reals $X$, the following are equivalent:
\be
\item $X$ satisfies $\ufin(\cO,\cO_\cF)$.
\item For each continuous image $Y$ of $X$ in $\NN$, there is $\sigma\in S_\N$ such that
$Y$ is bounded with respect to $\le_{\sigma\cF}$.
\ee
In particular, $\Bdd(\cF)$ implies $\ufin(\cO,\cO_\cF)$.
\end{thm}
\begin{proof}
$(2\Impl 1)$
Assume that $\cU_n$, $n\in\N$, are open covers of $X$,
which do not contain a finite subcover of $X$.
For each $n$, let $\tilde\cU_n$ be a countable cover refining $\cU_n$
such that all elements of $\tilde\cU_n$ are clopen and disjoint.
Enumerate $\tilde\cU_n$ bijectively as $\{C^n_m : m\in\N\}$.
Then the function $\Psi:X\to\NN$ defined by
$$\Psi(x)(n) = m\mbox{ such that }x\in C^n_m$$
is continuous, and therefore $\Psi[X]$ is bounded with respect to $\le_{\sigma\cF}$
for some $\sigma\in S_\N$.
Let $g\in\NN$ be a witness for that.
By induction on $n$, choose finite subsets $\cF_n\sbst\cU_n$ such that
$\Union_{m\le g(n)}C^n_m\sbst\bigcup\cF_n$, and such that
$\bigcup\cF_n$ is not equal to any $\bigcup\cF_k$ for $k<n$.\footnote{%
Since no $\cU_n$ contains a finite cover of $X$, we may achieve this as follows:
Choose a finite $\cA\sbst\cU_n$ such that $\Union_{m\le g(n)}C^n_m\sbst\bigcup\cA$.
For each $k<n$ take $x_k\in X\sm\bigcup\cF_k$.
Choose a finite $\cB\sbst\cU_n$ such that $\{x_1,\dots,x_{n-1}\}\sbst\bigcup\cB$,
and take $\cF_n=\cA\bigcup\cB$.
}
Consequently,
$$\{n : x\in \bigcup\cF_n\}\spst\{n : x\in \Union_{m\le g(n)}C^n_m\} = [\Psi(x)\le g]\in\sigma\cF$$
for each $x\in X$. Consequently, the (bijective!) enumeration $\{\bigcup\cF_{\sigma\inv(n)} : n\in\N\}$
witnesses that $\{\bigcup\cF_n : n\in\N\}$ is an $\cF$-cover.

$(1\Impl 2)$ Assume that $X$ satisfies $\ufin(\cO,\cO_\cF)$, and
let $Y$ be a continuous image of $X$ in $\NN$.
We may assume that each $f\in Y$ is increasing.
It is easy to see that the following holds.

\begin{lem}
$\ufin(\cO,\cO_\cF)$ is preserved under taking continuous images.\hfill\qed
\end{lem}

Thus, $Y$ satisfies $\ufin(\cO,\cO_\cF)$. Consider the open covers
$\cU_n=\{U^n_m : m\in\N\}$ of $Y$ defined by $U^n_m=\{f\in Y :
f(n)\le m\}$ (note that the elements $U^n_m$ are increasing with
$m$). There are two cases to consider.

Case 1: There is a strictly increasing sequence of natural
numbers $\seq{k_n}$, such that each $\cU_{k_n}$ contains
an element $U^{k_n}_{m_n}$ which is equal to $Y$.
Define $g(n)=m_n$ for each $n$.
Then for each $f\in Y$ and each $n$,
$f(n)\le f(k_n)\le m_n = g(n)$, that is, $[f\le g]=\N\in\cF$.

Case 2: There is $n_0$ such that for each $n\ge n_0$, $\cU_n$ does not contain $Y$ as an element.
Then by (1), there are finite subsets $\cF_n\sbst\cU_n$, $n\ge n_0$
such that $\cU=\{\bigcup\cF_n : n\ge n_0\}$ is an $\cF$-cover of $X$.
Let $h\in\NNup$ be such that $\{\bigcup\cF_{h(n)} : n\in\N\}$ is
a bijective enumeration of $\cU$.
Define $g(n)=\max\{m : U^{h(n)}_m\in\cF_{h(n)}\}$ for each $n$.
Then there is $\sigma\in S_\N$ such that
for each $f\in Y$, $\{n : f\in\bigcup\cF_{h(\sigma(n))}\}\in\cF$.
For each $n$ with $f\in\bigcup\cF_{h(\sigma(n))}$,
$$f(\sigma(n))\le f(h(\sigma(n)))\le g(\sigma(n)).$$
Thus $f\le_{\sigma\cF} g$ for all $f\in Y$.
\end{proof}

\begin{rem}\mbox{}
\be
\item By the methods of \cite{ZdImages}, Theorem \ref{bddimg} actually holds for arbitrary (not necessarily
zero-dimensional) subsets of $\R$.
\item One can characterize $\Bdd(\cF)$ by: For each sequence $\cU_n$ of open covers of $X$,
there exist finite subsets $\cF_n\sbst\cU_n$, $n\in\N$, such that for each $x\in X$,
$\{n : x\in\bigcup\cF_n\}\in\cF$.\hfill\qed
\ee
\end{rem}

Let $X$ be a set of reals.
In addition to $\gamma$-covers, the following type of covers plays a central role in the field:
An open cover $\cU$ of $X$ is an \emph{$\omega$-cover} of $X$ if $X$ is not in $\cU$ and for
each finite subset $F$ of $X$, there is $U\in\cU$ such that $F\subseteq U$.
Let $\Omega$ denote the collection of all countable open
$\omega$-covers of $X$.
For each \emph{filter} $\cF$, $\cO_\cF\sbst\Omega$.
Consequently, $\ufin(\cO,\cO_\cF)$ implies $\ufin(\cO,\Omega)$, which is strictly stronger
than Menger's property $\ufin(\cO,\cO)$ \cite{coc2}. In light of Theorem \ref{bddimg},
all of the examples shown to satisfy $\Bdd(\cF)$ for a filter $\cF$, satisfy $\ufin(\cO,\Omega)$.

\subsection{Finer distinction}
We now reveal the remainder of the framework of selection principles,
and apply the combinatorial approach to obtain a new result concerning these,
which further improves our earlier results.

\begin{figure}[!ht]
{\tiny
\begin{changemargin}{-3cm}{-3cm}
\begin{center}
$\xymatrix@C=7pt@R=20pt{
&
&
& \sr{\ufin(\cO,\Gamma)}{\fb~~ (18)}\ar[r]
& \sr{\ufin(\cO,\Tau)}{\max\{\fb,\s\}~~ (19)}\ar[rr]
&
& \sr{\ufin(\cO,\Omega)}{\fd~~ (20)}\ar[rrrr]
&
&
&
& \sr{\ufin(\cO,\cO)}{\fd~~ (21)}
\\
&
&
& \sr{\sfin(\Gamma,\Tau)}{\fbx{\fb}~~ (12)}\ar[rr]\ar[ur]
&
& \sr{\sfin(\Gamma,\Omega)}{\fd~~ (13)}\ar[ur]
\\
\sr{\sone(\Gamma,\Gamma)}{\fb~~ (0)}\ar[uurrr]\ar[rr]
&
& \sr{\sone(\Gamma,\Tau)}{\fbx{\fb}~~ (1)}\ar[ur]\ar[rr]
&
& \sr{\sone(\Gamma,\Omega)}{\fd~~ (2)}\ar[ur]\ar[rr]
&
& \sr{\sone(\Gamma,\cO)}{\fd~~ (3)}\ar[uurrrr]
\\
&
&
& \sr{\sfin(\Tau,\Tau)}{\fbx{\min\{\fb,\s\}}~~ (14)}\ar'[r][rr]\ar'[u][uu]
&
& \sr{\sfin(\Tau,\Omega)}{\fd~~ (15)}\ar'[u][uu]
\\
\sr{\sone(\Tau,\Gamma)}{\ft~~ (4)}\ar[rr]\ar[uu]
&
& \sr{\sone(\Tau,\Tau)}{\fbx{\ft}~~ (5)}\ar[uu]\ar[ur]\ar[rr]
&
& \sr{\sone(\Tau,\Omega)}{\fbox{\textbf{?}$(\fo)$}~~ (6)}\ar[uu]\ar[ur]\ar[rr]
&
& \sr{\sone(\Tau,\cO)}{\fbox{\textbf{?}$(\fo)$}~~ (7)}\ar[uu]
\\
&
&
& \sr{\sfin(\Omega,\Tau)}{\p~~ (16)}\ar'[u][uu]\ar'[r][rr]
&
& \sr{\sfin(\Omega,\Omega)}{\fd~~ (17)}\ar'[u][uu]
\\
\sr{\sone(\Omega,\Gamma)}{\p~~ (8)}\ar[uu]\ar[rr]
&
& \sr{\sone(\Omega,\Tau)}{\p~~ (9)}\ar[uu]\ar[ur]\ar[rr]
&
& \sr{\sone(\Omega,\Omega)}{\cov(\M)~~ (10)}\ar[uu]\ar[ur]\ar[rr]
&
& \sr{\sone(\cO,\cO)}{\cov(\M)~~ (11)}\ar[uu]
}$
\end{center}
\end{changemargin}
}

\caption{The surviving properties}\label{tauSch}
\end{figure}
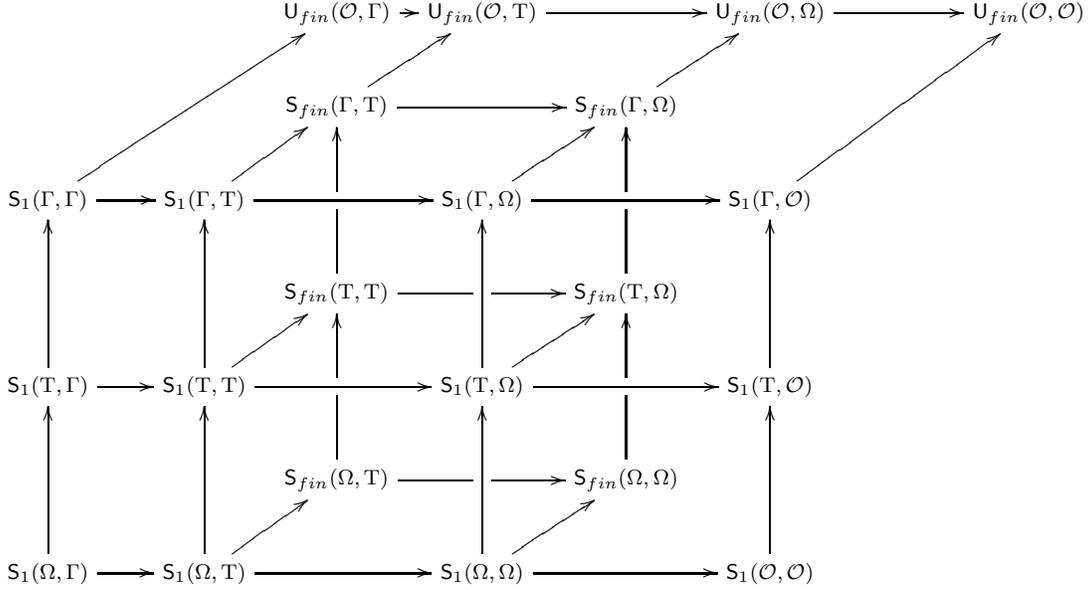

This framework was introduced by Scheepers in \cite{coc1, coc2} as a unified
generalization of several classical notions, and studied since in a long
series of papers by many mathematicians, see the surveys \cite{LecceSurvey, KocSurv, ict}.
Let $X$ be a set of reals.
An open cover $\cU$ of $X$ is a \emph{$\tau$-cover} of $X$ if each member of $X$ is covered
by infinitely many members of $\cU$,
and for each $x,y\in X$, at least one of the sets
$\{U\in\cU : x\in U\mbox{ and }y\nin U\}$ or
$\{U\in\cU : y\in U\mbox{ and }x\nin U\}$ is finite.
Let $\Tau$ denote the collection of all countable open
$\tau$-covers of $X$.
It is easy to see that
$$\Gamma\sbst\Tau\sbst\Omega\sbst\cO.$$
Let $\scrA$ and $\scrB$ be collections of covers of $X$.
In addition to $\ufin(\scrA,\scrB)$, we have the following selection hypotheses.
\bdesc
\item[$\sone(\scrA,\scrB)$]
For each sequence $\seq{\cU_n}$ of members of $\scrA$,
there exist members $U_n\in\cU_n$, $n\in\N$, such that $\setseq{U_n}\in\scrB$.
\item[$\sfin(\scrA,\scrB)$]
For each sequence $\seq{\cU_n}$
of members of $\scrA$, there exist finite (possibly empty)
subsets $\cF_n\sbst\cU_n$, $n\in\N$, such that $\Union_n\cF_n\in\scrB$.
\edesc
In addition to the Menger ($\ufin(\cO, \cO)$) and Hurewicz ($\ufin(\cO,\Gamma)$)
properties, several other properties of this form were studied in the past by
Rothberger ($\sone(\cO, \cO)$),
Arkhangel'ski\v{i} ($\sfin(\Omega,\Omega)$),\footnote{Arkhangel'ski\v{i} studied
``Menger property in all finite powers'', that was proved equivalent to
$\sfin(\Omega,\Omega)$ in \cite{coc2}.}
Gerlits and Nagy ($\sone(\Omega,\Gamma)$),
and Sakai ($\sone(\Omega,\Omega)$).
Many equivalences hold among these properties, and the surviving ones
appear in Figure \ref{tauSch} (where an arrow denotes implication) \cite{coc1, coc2, tautau}.

In \cite{coc2} it is proved that a set of reals $X$ satisfies $\sfin(\Omega,\Omega)$ if,
and only if, all finite powers of $X$ have the Menger property $\ufin(\cO,\cO)$.
By Corollary \ref{F+powers}, the examples involving filters (including those from Section \ref{algebra})
satisfy $\sfin(\Omega,\Omega)$. In Theorem \ref{UFscale} the example did not satisfy $\ufin(\cO,\Gamma)$.
We will improve that to find such an example which does not satisfy $\ufin(\cO,\Tau)$.
Since it is consistent that  $\ufin(\cO,\Gamma)$ is equivalent to $\ufin(\cO,\Tau)$ \cite{SF1}
and that $\ufin(\cO,\Omega)$ is equivalent to $\ufin(\cO,\cO)$ \cite{SF2}, our result is the best
possible with regards to Figure \ref{tauSch}.

Again, we will identify $\Inc$ with $\PN$.
We will use the following notion.
A family $Y\sbst\PN$ is \emph{splitting} if for each $A\in\roth$ there is $B\in Y$ such
that $A\cap B$ and $A\sm B$ are both infinite.
Recall that if $\fd$ is regular then there is an ultrafilter (necessarily nonmeager) $\cF$
satisfying $\bofF=\fd$.

\begin{thm}\label{Tauscale}
Assume that $\fd$ is regular.
Then for each nonmeager filter $\cF$ with $\bofF=\fd$,
there is a cofinal $\bofF$-scale $S=\{a_\alpha : \alpha<\fd\}\sbst\ici$ such that:
\be
\item All finite powers of the set $X=S\cup Q$ satisfy $\Bdd(\cF)$,
but
\item The homeomorphic copy $\tilde X=\{x\comp : x\in X\}$ of $X$
is a splitting and unbounded (with respect to $\le^*$) subset of $\roth$.
\ee
\end{thm}
\begin{proof}
For each $h\in\NNup$, let
$$\cA_h=\left\{\Union_{n\in A}\intvl{h}{n} : A\in\ici\right\}.$$
$\cA_h$ is homeomorphic $\ici$, and is therefore analytic.

\begin{lem}\label{Abdd}
For each $h\in\NNup$ and each $f\in\NNup$, there is $a\in\cA_h$ such that $f\le_\cF a$.
\end{lem}
\begin{proof}
Clearly, $\cA_h$ is not $\le^*$-bounded.
Apply Lemma \ref{nonM1}.
\end{proof}

\begin{lem}\label{AKbdd}
For all $h,f,g\in\NNup$,
there is $a\in\cA_h$ such that $f\le_\cF a$ and $a\comp\not\le^* g$.
\end{lem}
\begin{proof}
Let $q\in\NNup$ be such that for each $a\in\ici$ with
$a\le^* g$,
$a$ intersects all but finitely many of the intervals $\intvl{q}{n}$.
(E.g., define inductively $q(0) = g(0)$ and $q(n+1) = g(q(n))+1$.)
We may assume that $\im q\sbst \im h$, and therefore $\cA_q\sbst\cA_h$.
By Lemma \ref{Abdd}, there is $a\in\cA_q$ such that $f\le_\cF a$.
Since $a\comp\in\cA_q$, it misses infinitely many intervals $\intvl{q}{n}$,
and therefore $a\comp\not\le^* g$.
\end{proof}

Let $\{d_\alpha : \alpha<\fd\}\sbst\NNup$ be such that for each $A\in\roth$
there is $\alpha<\fd$ such that $|A\cap \intvl{d_\alpha}{n}|\ge 2$
for all but finitely many $n$ \cite{BlassHBK}.
In particular, $\{d_\alpha : \alpha<\fd\}$ is dominating.

For each $\alpha<\bofF=\fd$ inductively, do the following:
Choose $f\in\NNup$ which is a $\le_\cF$-bound of $\{a_\beta : \beta<\alpha\}$.
Use Lemma \ref{AKbdd} to choose $a_\alpha\in\cA_{d_\alpha}$ such that
$\max\{f,d_\alpha\}\le_\cF a_\alpha$ and $a_\alpha\comp\not\le^* d_\alpha$.
Since $\cF$ is a filter, $a_\beta\le_\cF a_\alpha$ for each $\beta<\alpha$.

$S=\{a_\alpha : \alpha<\fd\}$ is a cofinal $\bofF$-scale, and thus by Theorem
\ref{cofscale2}, all finite powers of $X=S\cup Q$ satisfy $\Bdd(\cF)$.
As for each $\alpha<\fd$ we have $a_\alpha\comp\not\le^* d_\alpha$,
$\tilde X$ is unbounded.
To see that it is splitting, let $b\in\roth$ and choose $\alpha$ such that
$b$ intersects $\intvl{d_\alpha}{n}$ for all but finitely $n$.
Since $a_\alpha\in\cA_{d_\alpha}$, $a_\alpha$ splits $b$.
\end{proof}

According to \cite{tautau}, a subset $Y$ of $\NN$
has the \emph{excluded-middle property}
if there exists $g\in\NN$ such that:
\be
\item for each $f\in Y$, the set $[f<g]$ is infinite; and
\item for all $f,h\in Y$ at least one of the sets
$[f<g\le h]$ and $[h<g\le f]$ is finite.
\ee
In Theorem 3.11 and Remark 3.12 of \cite{tautau} it is
proved that if $Y$ satisfies $\ufin(\cO,\Tau)$, then
all continuous images of $Y$ in $\NN$ have
the excluded-middle property.

\begin{cor}\label{TauEx}
Assume that $\fd$ is regular.
Then for each nonmeager filter $\cF$ with $\bofF=\fd$,
there is a set of reals $Y\sbst\NN$ such that:
\be
\item All finite powers of $Y$ satisfy $\Bdd(\cF)$, but
\item $Y$ does not have the excluded-middle property. In particular,
$Y$ does not satisfy $\ufin(\cO,\Tau)$.
\ee
\end{cor}
\begin{proof}
Let $\tilde X$ be the set from Theorem \ref{Tauscale}(2).
Define continuous functions $\Psi_\ell:\tilde X^2\to\NN$, $\ell=0,1$, by
$$\Psi_0(x,y)(n) = \begin{cases}
x(n) & n\in y\\
0 & n\nin y
\end{cases}
;\qquad
\Psi_1(x,y)(n) = \begin{cases}
0 & n\in y\\
x(n) & n\nin y
\end{cases}
$$
Take $Y=\Psi_0[\tilde X^2]\cup\Psi_1[\tilde X^2]$.
Each finite power of $Y$ is a finite union of continuous images
of finite powers of $\tilde X$. Consequently, all finite powers of $Y$
satisfy $\Bdd(\cF)$.

The argument in the proof of Theorem 9 of \cite{ShTb768} shows that
$Y$ does not have the excluded middle property.
\end{proof}

We obtain the following.

\begin{thm}
There exists a set of reals $X$ satisfying $\sfin(\Omega,\Omega)$ but
not $\ufin(\cO,\Tau)$.
\end{thm}
\begin{proof}
The proof is dichotomic.
If $\cf(\fd)=\fd$, use Corollary \ref{TauEx}.
Otherwise, $\cf(\fd)<\fd$.
As $\max\{\fb,\s\}\le\cf(\fd)$ ($\s\le\cf(\fd)$ is proved in \cite{Mild01}),
$\max\{\fb,\s\}<\fd$.
As the critical cardinalities of $\ufin(\cO,\Tau)$ and $\sfin(\Omega,\Omega)$
are $\max\{\fb,\s\}$ \cite{ShTb768} and $\fd$ \cite{coc2}, respectively,
we can take $Y$ to be any witness for the first of these two assertions.
\end{proof}

By the arguments of Section \ref{algebra}, we have the following.
\begin{cor}\label{strongfields}
Assume that $\mathbb{F}$ is a nondiscrete, separable,
completely metrizable field, and $\Q$ is a countable dense subfield of $\mathbb{F}$.
\be
\item If $\fd$ is regular, then for each nonmeager filter $\cF$ with $\bofF=\fd$,
there is $X\sbst\mathbb{F}$ such that:
\be
\item All finite powers of $\Q(X)$ satisfy $\Bdd(\cF)$, but
\item $\Q(X)$ does not satisfy $\ufin(\cO,\Tau)$.
\ee
\item There exists $X\sbst\mathbb{F}$ such that $\Q(X)$ satisfies $\sfin(\Omega,\Omega)$ but
not $\ufin(\cO,\Tau)$.\hfill\qed
\ee
\end{cor}

Readers not familiar with forcing can safely skip the following remark.

\begin{rem}
The constructions in this section can be viewed as an extraction
of the essential part in the forcing-theoretic construction obtained by
adding $\fc$ many Cohen reals to a model of set theory, and letting
$X$ be the set of the added Cohen reals. Since Cohen reals
are not dominating, all finite powers of $X$ will have Menger's property.
It is also easy to see that $X$ will not satisfy the excluded-middle property,
e.g., using the reasoning in \cite{ShTb768}.
See \cite{JORG} for these types of constructions, but note that they only yield
consistency results.
\end{rem}


\section{Towards semifilters again}\label{SFA}

We strengthen the solution to the Hurewicz Problem as follows.

\begin{thm}
Assume that $P$ is a nowhere locally compact Polish space,
and $\cS$ is a nonmeager bounding \emph{semifilter} such that
$\bof(\cS)=\fd$.
Then there is a subspace $X$ of $P$ such that:
\be
\item All finite powers of $X$ have Menger's property,
\item $X$ satisfies $\Bdd(\cS)$; and
\item $X$ does not have the Hurewicz property.
\ee
\end{thm}
\begin{proof}
As pointed out in Remark \ref{Polish}, it suffices to consider the
case $P=\ici\cup\Fin$, in a disguise of our choice.
We give an explicit construction in the case that $\fd$ is regular.
The remaining case, being ``rare'' but consistent, is trivial.

A family $\cF\sbst\roth$ is \emph{centered} if each finite subset
of $\cF$ has an infinite intersection.
Centered families generate filters by taking finite intersections and
closing upwards. We will denote the generated filter by $\langle \cF\rangle$.
For $Y\sbst\NN$, let $\maxfin Y$ denote its closure under
pointwise maxima of finite subsets.

We construct, by induction on $\alpha<\fd$, a filter $\cF$ with $\bofF=\fd$
and a $\bofF$-scale $\{a_\alpha : \alpha<\fd\}\sbst\ici$
which is also a cofinal $\bof(\cS)$-scale.

Let $\{d_\alpha : \alpha<\fd\}\sbst\NN$ be dominating, and assume
that $a_\beta$ are defined for each $\beta<\alpha$.
Let
\begin{eqnarray*}
\cA_\alpha & = & \maxfin\{d_\beta, a_\beta : \beta<\alpha\},\\
\tilde\cF_\alpha & = & \Union_{\beta<\alpha}\cF_\beta;\\
\cG_\alpha & = & \{f\circ b : f\in\cA_\alpha, b\in\tilde\cF_\alpha\}.
\end{eqnarray*}
We inductively assume that $\cF_\beta$, $\beta<\alpha$, is an increasing
chain of filters such that $|\cF_\beta|\le|\beta|$ for each $\beta<\alpha$.
This implies that $|\cG_\alpha|\le|\alpha|<\fd$.

As $\cS$ is nonmeager-bounding, there exists
a $\le_\cS$-bound $a_\alpha$ of $\cG_\alpha$
such that $a_\alpha\comp\not\le^* d_\alpha$.
Define
$$\cF_\alpha = \langle \tilde\cF_\alpha\cup\{[f\le a_\alpha] : f\in\cA_\alpha\}\rangle.$$
We must show that $\cF_\alpha$ remains a filter.
First, assume that there are $b\in\tilde\cF_\alpha$ and $f\in\cA_\alpha$
such that $b\cap[f\le a_\alpha]$ is finite.
Then $a_\alpha\le a_\alpha\circ b <^* f\circ b\in\cG_\alpha$,
a contradiction.
Now, for each $b\in\cF_\alpha$ and $f_1,\dots f_k\in\cA_\alpha$,
we have that $f=\max\{f_1,\dots,f_k\}\in\cA_\alpha$, and therefore
$$b\cap\bigcap_{i=1}^k [f_i\le a_\alpha]= b\cap [f\le a_\alpha]$$
is infinite.

Take $S=\{a_\alpha : \alpha<\fd\}$, and
$\cF=\Union_{\alpha<\fd}\cF_\alpha$.
By the construction, $S$ is a cofinal $\bof(\cS)$-scale.
By Theorem \ref{cofscale1},
$X=S\cup Q$ satisfies $\Bdd(\cS)$.
For all $\alpha<\beta<\fd$, $a_\alpha\le_\cF a_\beta$.
We claim that if $\fd$ is regular, then $\bof(\cF)=\fd$.
Indeed, assume that $Y\sbst\NN$ has cardinality less than $\fd$.
As $\fd$ is regular, there exists $\alpha<\fd$ such that
each $f\in Y$ is $\le^*$-bounded by some $d_\beta$, $\beta<\alpha$.
As $a_\alpha$ is a $\le_\cF$-bound of $\{d_\beta : \beta<\alpha\}$,
it is a $\le_\cF$-bound of $Y$.
We get that $S$ is also a cofinal $\bofF$-scale.
By Theorem \ref{cofscale2},
all finite powers of $X$ satisfy $\Bdd(\cF)$
(and, in particular, Menger's property).

Finally, since $\{x\comp : x\in X\}$ is an unbounded subset of $\NN$, $X$ does not have
the Hurewicz property.
\end{proof}

\section{Topological Ramsey theory}\label{Ramsey}

Most of our constructions can be viewed as examples in topological Ramsey theory.
We explain this briefly.
The following partition relation, motivated
by a study of Baumgartner and Taylor in Ramsey theory \cite{BT78},
was introduced by Scheepers in \cite{coc1}:
\bdesc
\item[$\scrA\rightarrow \lceil\scrB\rceil^2_k$]
For each $\cU\in\scrA$ and each $f:[\cU]^2 \rightarrow \{0,\dots,k-1\}$,
there exist $\cV\sbst \cU$ such that $\cV\in\scrB$, $j\in\{0,\cdots,k-1\}$, and a partition $\cV = \bigcup_n\cF_n$
of $\cV$ into finite sets, such that for each $\{A,B\}\in [\cV]^2$
such that $A$ and $B$ are not from the same $\cF_n$, $f(\{A,B\})=j$.
\edesc

Menger's property is equivalent to
$(\forall k)\ \Omega \rightarrow \lceil\cO\rceil^2_k$ \cite{OpPar}, and
having the Menger property
in all finite powers is equivalent to
$(\forall k)\ \Omega\rightarrow\lceil\Omega\rceil^2_k$ \cite{coc3}.

A cover $\cU$ of $X$ which does not contain a finite subcover
is \emph{$\gamma$-groupable} if
there is a partition of $\cU$ into finite sets, $\cU = \Union_n\cF_n$, such  that
$\{\cup\cF_n : n\in\N\}$ is a $\gamma$-cover of $X$.
Denote the collection of $\gamma$-groupable open covers of $X$ by $G(\Gamma)$.

The Hurewicz property is equivalent to
$(\forall k)\ \Omega\rightarrow\lceil G(\Gamma)\rceil^2_k$,
and having the Hurewicz property
in all finite powers is equivalent to
$(\forall k)\ \Omega\rightarrow\lceil\Omega^{\grbl}\rceil^2_k$,
where $\Omega^{\grbl}$ denotes covers with partition $\Union_n\cF_n$
into finite sets such that for each finite $F\sbst X$ and all
but finitely many $n$, there is $U\in\cF_n$ such that $F\sbst U$ \cite{coc7, QRT, RTG}.

Clearly, $\Omega^{\grbl}\sbst\Omega\cap G(\Gamma)$.

We state only three of our results using this language, leaving the
statement of the remaining ones to the reader.

\begin{thm}
\mbox{}
\be
\item The sets $X$ constructed in Theorem \ref{answerHure} satisfy
$(\forall k)\ \Omega \rightarrow \lceil\cO\rceil^2_k$
but not
$(\forall k)\ \Omega\rightarrow\lceil G(\Gamma)\rceil^2_k$.
\item The fields $\Q(X)$ constructed in Theorem \ref{PlewPol}
satisfy $(\forall k)\ \Omega\rightarrow\lceil\Omega^{\grbl}\rceil^2_k$
but are not $\sigma$-compact.
\item The fields $\Q(X)$ constructed in Theorem \ref{strongfields}
satisfy $(\forall k)\ \Omega\rightarrow\lceil\Omega\rceil^2_k$ but not
$(\forall k)\ \Omega\rightarrow\lceil\Omega^{\grbl}\rceil^2_k$
(or even $(\forall k)\ \Omega\rightarrow\lceil G(\Gamma)\rceil^2_k$).
\ee
\end{thm}

\subsection{Strong measure zero and Rothberger fields}
We need not stop at the decidable case.
According to Borel \cite{Borel},
a set of reals $X$ has \emph{strong measure zero} if
for each sequence of positive reals $\seq{\epsilon_n}$,
there exists a cover
$\setseq{I_n}$ of $X$ such that for each $n$, the diameter of $I_n$ is smaller than $\epsilon_n$.
This is a very strong property, and Borel conjectured that every strong measure zero set of reals
is countable. This was proved consistent by Laver \cite{Laver}.

Rothberger's property $\sone(\cO,\cO)$ implies strong measure zero, and its critical cardinality
is $\cov(\cM)$, the minimal cardinality of a cover of the real line by meager sets.
By known combinatorial characterizations \cite{barju},
if $\fb$ is not greater than the minimal cardinality of a set of reals
which is not of strong measure zero, then $\fb\le\cov(\cM)$.
In the following theorem, any embedding of $\Inc$ in $\R$
can be used.
\begin{thm}
If $\fb\le\cov(\cM)$,
then the fields $\Q(X)$ constructed in Theorem \ref{PlewPol}
satisfy $\sone(\Omega,\Omega^{\grbl})$ (equivalently,
all finite powers of $\Q(X)$ satisfy the Hurewicz property as well as
Rothberger's property \cite{coc7}).
\end{thm}
\begin{proof}
Since $X=S\cup Q$ is $\fb$-concentrated on the countable set $Q$, it
satisfies---by the assumption on $\fb$---Rothberger's property $\sone(\cO,\cO)$.
As all finite powers of $X$ have the Hurewicz property,
we have by Theorem 4.3 of \cite{prods} that $X$ satisfies $\sone(\Omega,\Omega^{\grbl})$.
By the arguments in the proof of Theorem \ref{basicfields}(1),
all finite powers of $\Q(X)$ satisfy $\sone(\cO,\cO)$.
Theorem \ref{basicfields}(2) tells us that all finite powers
of $\Q(X)$ also satisfy the Hurewicz property, so we are done.
\end{proof}

The following partition relation \cite{coc1} is a natural extension of Ramsey's.
\bdesc
\item[$\scrA\to(\scrB)^n_k$]
For each $\cU\in\scrA$ and $f:[\cU]^n\to\{1,\ldots,k\}$,
there exist $j$ and $\cV\sbst\cU$ such that $\cV\in\scrB$ and $f\rest[\cV]^n\equiv j$.
\edesc
Using this notation, Ramsey's Theorem is $(\forall n,k)\ \roth\to(\roth)^n_k$.

In \cite{coc7} it is proved that $\sone(\Omega,\Omega^{\grbl})$ is equivalent to
$(\forall n,k)\ \Omega\rightarrow(\Omega^{\grbl})^n_k$.

\begin{cor}
If $\fb\le\cov(\cM)$,
then the fields $\Q(X)$ constructed in Theorem \ref{PlewPol}
satisfies $(\forall n,k)\ \Omega\rightarrow(\Omega^{\grbl})^n_k$.
\end{cor}

\section{Some concluding remarks}

Using filters in the constructions allowed avoiding some of the technical
aspects of earlier constructions and naturally obtain examples for the Menger and Hurewicz
Conjectures which possess an algebraic structure.
The extension to semifilters is essential for the consideration of the Menger
and Hurewicz properties in terms of boundedness on ``large'' sets of natural numbers.
While making some of the proofs more difficult, it
seems to have provided the natural solution of the Hurewicz Problem, and allowed
its strengthening in several manners.

Chaber and Pol asked us about the difference in strength between the construction
in \cite{ideals} (corresponding to item (2) in Corollary \ref{BFpowers})
and their dichotomic construction \cite{ChaPol} (Theorem \ref{CPThm}).
The answer is now clear:
The set from \cite{ideals} has the Hurewicz property, and the (dichotomic) set from
\cite{ChaPol} has the Menger property but not the Hurewicz property.

Previous constructions (dichotomic or ones using additional hypotheses) which made various
assumptions on the cardinal $\fd$ can now be viewed as a ``projection'' of the constructions
which only assume that $\fd$ is regular. While giving rise (in a dichotomic manner) to
ZFC theorems, the possibility to eliminate the dichotomy
in Theorem \ref{UFscale} and its consequences
without making any additional hypotheses remains open.
It may be impossible.

\bigskip\subsection*{Acknowledgements}
The proof of Theorem \ref{PlewPol} evolved from a suggestion of Roman Pol,
which we gratefully acknowledge.
The proof of this theorem also uses a result of Plewik, whom we acknowledge for
letting us know about it.
We also thank Taras Banakh and Micha\l{} Machura for reading the paper and
making comments.
A special thanks is owed to Marion Scheepers for his detailed comments
on this paper.

\ed